\documentclass{elsart}

\usepackage{amssymb}
\usepackage[french]{babel}

\newtheorem{theoreme}{Th\'eor\`eme}[section]
\newtheorem{defin}[theoreme]{D\'efinition}
\newtheorem{propo}[theoreme]{Proposition}

\newtheorem{corol}[theoreme]{Corollaire}
\newtheorem{lemmeb}{Lemme.}

\newenvironment{demo}{\noindent {\sl Preuve}. \ }{\qed}

\newtheorem{stheoreme}{Th\'eor\`eme }[subsection]
\newtheorem{sdefin}[stheoreme]{D\'efinition }
\newtheorem{spropo}[stheoreme]{Proposition}
\newtheorem{slemme}[stheoreme]{Lemme}
\newtheorem{scorol}[stheoreme]{Corollaire}
\newtheorem{propb}{Proposition.}

\newtheorem{definb}{D\'efinition.}

\newtheorem{corolb}{Corollaire.}

\usepackage{psfig}

\def\winf{\mathop{\hbox{\textrm w-infdeg}}\nolimits}
\def\wsup{\mathop{\hbox{\textrm w-supdeg}}\nolimits}
\def\winfpr{\mathop{\hbox{\textrm{\emph{w-infdeg}}}}\nolimits}
\def\wsuppr{\mathop{\hbox{\textrm{\emph{w-supdeg}}}}\nolimits}

\def\soul{\underline}
\def\tl{{\mathcal O}_{\L}}
\def\AG{{\mathcal A}(\G)}
\def\Alr{S_{n,r}^{\L}}

\def\smmath{\tiny\textrm}
\def\pprod{\prod\limits}
\def\ssum{\sum\limits}

\font\tenCal=cmsy10
\newfam\Calfam
\textfont\Calfam=\tenCal

\def\soul{\underline}

\def\AG{{\mathcal A}(\G)}
\def\Alr{S_{n,r}^{\L}(k)}

\def\Frac{{\mathcal F}{\textrm{rac}}}

\def\Aut{\mathop{\hbox{\textrm Aut}}\nolimits}
\def\End{\mathop{\hbox{\textrm End}}\nolimits}

\def\qed{\hfill{$\sqcap\!\!\!\!\sqcup$}}

\def\wt{\widetilde}
\def\what{\widehat}
\def\oif{\Rightarrow}
\def\surl{\overline}
\def\gkd{\mathop{\hbox{\textrm GKdim}}\nolimits}
\def\gktd{\mathop{\hbox{\textrm GKtrdeg}}\nolimits}
\def\gkdpr{\mathop{\hbox{\textrm{\emph{GKdim}}}}\nolimits}
\def\gktdpr{\mathop{\hbox{\textrm{\emph{GKtrdeg}}}}\nolimits}

\def\An{A_n^{\bar{q} ,\Lambda}}

\def\Bn{B_n^{\bar{q} ,\Lambda}}

\def\Dr{{\mathcal F}{\rm rac}(S^{\L}_{n,r}(k))}

\def\qb{\bar{q}}
\def\L{\Lambda}
\def\l{\lambda}
\def\a{\alpha}
\def\b{\beta}

\def\G{\Gamma}
\def\d{\delta}

\def\blanc{\ \ \ \ \ \ \ \ \ \ }
\def\esp{\ \ \ \ \ }

\begin{document}

\begin{frontmatter}

\title{Equivalence rationnelle d'alg\`ebres polynomiales classiques et quantiques}
\author{Lionel Richard
}
\address{Institut Girard Desargues,
Universit\'e Claude Bernard,  69622 Villeurbanne Cedex, France}
\ead{Lionel.Richard@igd.univ-lyon1.fr}
\journal{Journal of Algebra}
\date{7 Juillet 2003}

\begin{abstract}
Cet article est consacr\'e \`a l'\'etude de l'\'equivalence rationnelle des alg\`ebres non commutatives de polyn\^omes dans un cadre englobant \`a la fois le probl\`eme classique de Gelfand-Kirillov et son analogue quantique.
On introduit dans ce contexte mixte
 une classe d'alg\`ebres de r\'ef\'erence et on d\'efinit deux nouveaux invariants permettant  de s\'eparer les corps de fractions de ces alg\`ebres \`a isomorphisme pr\`es.
L'un est li\'e \`a la notion de sous-tore quantique maximal simple, l'autre est un invariant dimensionnel mesurant {\it via} certains plongements le caract\`ere classique, en terme d'alg\`ebre de Weyl, des corps consid\'er\'es.
 A titre d'application on en d\'eduit des r\'esultats concernant  l'\'equivalence rationnelle des alg\`ebres de Weyl quantiques multiparam\'etr\'ees.\\[.2cm]
{\parindent=0pt \bf Abstract}\\[.2cm]
{\parindent=0pt This article is devoted to rational equivalence for non-commutative polynomial algebras in a context including both  the classical Gelfand-Kirillov problem and its quantum version.
We introduce in this ``mixed'' context  some reference algebras and define two new invariants allowing us to separate the fraction fields of these algebras up to isomorphism.
The first one is linked to the notion of maximal simple quantum sub-torus, and the second one is a dimensionnal invariant measuring the classical character (in terms of Weyl algebras) of the skew-fields into consideration.
As an application we obtain results concerning  the rational equivalence of multiparametrized quantum Weyl algebras.}
\end{abstract}



\begin{keyword} Rational equivalence, iterated Ore extensions, classical and quantum Weyl algebras
\MSC 16K40, 16W35, 16S32
\end{keyword}

\end{frontmatter}

\section*{Introduction.}

L'\'etude des corps de fonctions rationnelles non commutatifs
(en un nombre fini de variables) reste fortement impr\'egn\'ee des
probl\'ematiques de l'article fondateur \cite{GK} de I.M. Gelfand et A.A. Kirillov
de 1966.
Deux d'entre elles retiendront notre attention dans ce qui suit.
\begin{itemize}
\item
La premi\`ere concerne la question de d\'eterminer des classes d'alg\`ebres non commutatives de polyn\^omes
(noeth\'eriennes, int\`egres) rationnellement \'e\-qui\-va\-len\-tes \`a
certains types canoniques.
\item
La seconde consiste \`a d\'egager des invariants permettant la classification
ou la s\'eparation de ces alg\`ebres canoniques entre elles \`a \'equivalence
rationnelle pr\`es.
\end{itemize}
Dans le contexte ``classique", le premier probl\`eme s'articule dans sa
formulation originelle autour de l'\'equivalence rationnelle des alg\`ebres
enveloppantes d'alg\`ebres de Lie alg\'ebriques avec les alg\`ebres de Weyl
sur un anneau de polyn\^omes. 
Un contre-exemple a \'et\'e produit dans l'article \cite{AOVB}, o\`u l'on trouvera aussi des r\'ef\'erences exhaustives sur les nombreuses situations o\`u la conjecture est v\'erifi\'ee.
Le second probl\`eme est enti\`erement r\'esolu par deux invariants dimensionnels
entiers : le degr\'e de transcendance du centre et la notion de degr\'e
de transcendance non-commutatif de Gelfand et Kirillov, qui suffisent \`a
s\'eparer \`a isomorphisme pr\`es les corps de fractions des alg\`ebres
de Weyl (corps de Weyl).
L'ensemble constitue ce qu'il est d'usage de d\'esigner sous le nom
de probl\`eme de Gelfand et Kirillov.

Dans le contexte ``quantique" le premier probl\`eme (formul\'e  par exemple en \cite{AD1},  \cite{IM}, \cite{Josep}, \cite{Cal}, \cite{Jslqwa}, etc.) a conduit au cours
des ann\'ees 1990 \`a prouver l'\'equivalence rationnelle de vastes classes
d'alg\`ebres ``quantiques"  avec des espaces
affines quantiques.
La litt\'erature concernant la version quantique de la conjecture de Gelfand-Kirillov
     est abondante (pour une bibliographie plus compl\`ete nous renvoyons le lecteur \`a l'article \cite{AD3} et au paragraphe II.10.4 de l'ouvrage \cite{BG}), et il semble s'agir d'une propri\'et\'e tout \`a
     fait g\'en\'erale concernant de vastes classes d'alg\`ebres
      noe\-th\'e\-rien\-nes int\`egres (\cite{Pa2}, \cite{Ca2}).
 Le second probl\`eme est moins simple que dans le cas classique, les invariants
dimensionnels pr\'ec\'edents ne suffisant pas \`a s\'eparer \`a isomorphisme
pr\`es les corps de fractions des espaces affines quantiques, du
fait de leur param\'etrisation  par toute une matrice de
scalaires non-nuls. Sur ce point, mentionnons les r\'esultats obtenus dans le cas uniparam\'etr\'e par A.N. Panov (\cite{Pa1}), dans le cas ``g\'en\'erique'' par V.A. Artamonov (\cite{A}, \cite{A2}), et enfin le th\'eor\`eme 4.2 de \cite{Retq},  rappel\'e ici en \ref{pigne}, portant sur l'\'equivalence rationnelle des tores quantiques multiparam\'etr\'es (alg\`ebres de polyn\^omes de Laurent non-commutatifs) simples.

L'article \cite{CGG} de L.J. Corwin, I.M. Gelfand et R. Goodman semble \^etre l'un des
premiers \`a poser explicitement des questions relatives \`a l'\'equivalence rationnelle dans des situations mixtes (c'est-\`a-dire relevant simultan\'ement des contextes classique et quantique),
 en prenant pour alg\`ebres
de r\'ef\'erence des produits tensoriels d'alg\`ebres de Weyl et de plans
quantiques. Mais c'est dans l'article \cite{AD3} de J. Alev et F. Dumas que
sont obtenus sur ces questions des r\'esultats significatifs avec,
d'une part la classification \`a \'equivalence rationnelle pr\`es de ces
alg\`ebres de r\'ef\'erence, et d'autre part la production de contre-exemples
d'alg\`ebres de polyn\^omes \`a la fois classiques et quantiques mais non
rationnellement \'equivalentes \`a un produit tensoriel d'une alg\`ebre de
Weyl et d'un espace affine quantique. A la fois par certains outils qui y
sont introduits et par certaines questions ouvertes qui y sont formul\'ees, l'article \cite{AD3} peut \^etre consid\'er\'e comme un des points de d\'epart du pr\'esent travail.

\medskip

Dans la premi\`ere section de cet article nous introduisons les alg\`ebres polynomiales mixtes crois\'ees $S_{n,r}^{\L}(k)$ (d\'ej\`a mentionn\'ees dans \cite{Rhhqw}), extensions it\'er\'ees de Ore param\'etr\'ees par deux entiers $n$ et $r$ et une matrice $\L\in M_n(k^*)$ en en donnant la d\'efinition, des exemples et les premi\`eres propri\'et\'es.
Le r\^ole de ``standard'' de ces alg\`ebres pour  les situations \`a la fois classiques et quantiques est mis en \'evidence \`a la section 2. On y montre en effet 
que 
 toute
 alg\`ebre de polyn\^omes en $N$
variables $x_1,\ldots,x_N$ telle que la relation de commutation
entre deux g\'en\'erateurs ne soit que de l'un des deux types
suivants, quantique ($x_ix_j=\lambda_{i,j}x_jx_i$ avec $\lambda_{i,j}\in k^*$)
ou classique ($x_ix_j-x_jx_i=1$), est n\'ecessairement isomorphe \`a une alg\`ebre polynomiale mixte crois\'ee. La preuve se fait par une m\'ethode combinatoire sur les
graphes associ\'es \`a de tels syst\`emes de g\'en\'erateurs.
Via des m\'ethodes de 
plongements dans des corps de s\'eries de Malcev-Neumann (m\'ethodes maintenant
 standard), 
on calcule \`a la section 3 pour les corps de fractions de ces alg\`ebres
leur centre, ainsi que  les invariants rationnels $E$ et $G$ introduits dans \cite{AD1}, permettant de s\'eparer les situations purement classiques des situations purement quantiques.
L'\'etude des situations mixtes conduit \`a introduire deux nouveaux invariants.  On montre dans la section 4 que
le tore quantique param\'etr\'e par la matrice $\Lambda$, lorsqu'il est simple, est un invariant de ces corps, appel\'e  sous-tore quantique
maximal simple.
La section 5 est consacr\'ee \`a la d\'efinition et au calcul d'un autre invariant, dimensionnel celui-ci  : le w-degr\'e sup\'erieur, \'egal \`a la dimension de la plus petite alg\`ebre de Weyl n\'ecessaire au plongement d'une alg\`ebre donn\'ee dans le corps de fractions du produit tensoriel de cette alg\`ebre de Weyl avec un espace affine quantique. 
On montre alors que ce w-degr\'e sup\'erieur n'est autre pour le corps de fractions d'une alg\`ebre $S_{n,r}^{\Lambda}(k)$ que le double de l'entier $r$ param\'etrant cette alg\`ebre.
Sous la seule hypoth\`ese (raisonnable au
vu du th\'eor\`eme \ref{thpr}) de la simplicit\'e du tore quantique param\'etr\'e par $\Lambda$, on en d\'eduit \`a la section 6 un syst\`eme
de conditions n\'ecessaires  pour l'isomorphisme de  corps de fonctions rationnelles mixtes crois\'es.
Ces conditions s'av\`erent \^etre \'egalement suffisantes dans le cas dit semi-classique (o\`u $n=r$),
correspondant aux corps de fractions des alg\`ebres d'op\'erateurs diff\'erentiels eul\'eriens sur l'espace quantique param\'etr\'e par $\L$, redonnant ainsi un r\'esultat d\'ej\`a d\'emontr\'e dans \cite{Retq}.
Enfin la derni\`ere section est consacr\'ee \`a l'application des r\'esultats pr\'ec\'edents \`a  une classe d'alg\`ebres significative admettant des localisations communes avec des alg\`ebres polynomiales mixtes crois\'ees : les alg\`ebres de Weyl quantiques multiparam\'etr\'ees.

Dans toute la suite $k$ d\'esigne un corps de caract\'eristique z\'ero, et tous les morphismes consid\'er\'es sont des morphismes de $k$-alg\`ebres.

\medskip

{\parskip=0pt
\tableofcontents
}


\section{Alg\`ebres polynomiales mixtes crois\'ees.} \label{defsnrl}
 
\subsection{D\'efinitions et premi\`eres propri\'et\'es.} \label{defqw}

Afin de fixer les notations nous rappelons les d\'efinitions suivantes (voir par exemple le  chapitre I.2 de \cite{BG}).

\begin{sdefin} \label{dq}
\begin{enumerate}
\item Soit $n\geq 1$ un entier.
Une matrice $\L=(\l_{i,j})\in M_n(k^*)$ est dite multiplicativement
antisym\'etrique si ses coefficients v\'e\-ri\-fient 
$$\l_{i,j}\l_{j,i}=\l_{i,i}=1 \ \textrm{ pour tous } 1\leq i,j\leq n.$$
\item Soit $\L\in M_n(k^*)$ une matrice multiplicativement antisym\'etrique.
On appelle espace affine quantique param\'etr\'e par $\L$, et on note 
$${\mathcal O}_{\L}(k^n)=k_{\L}[y_1,\ldots,y_n],$$
 la $k$-alg\`ebre engendr\'ee par $n$ g\'en\'erateurs $y_1,\ldots,y_n$ avec les relations 
$$y_iy_j=\l_{i,j}y_jy_i \ \textrm{ pour tous } i,j.$$
\item  On appelle tore quantique  param\'etr\'e par $\L$, et on note 
$${\mathcal O}_{\L}({k^*}^n)=k_{\L}[y_1^{\pm 1},\ldots,y_n^{\pm 1}],$$
la localisation de l'espace affine quantique  $k_{\L}[y_1,\ldots,y_n]$ en la partie multiplicative engendr\'ee par les \'el\'ements  normaux $y_1,\ldots,y_n$.
\item On appelle corps de fonctions rationnelles quantique param\'etr\'e par $\L$ le corps non commutatif
$$k_{\L}(y_1,\ldots,y_n)=\Frac({\mathcal O}_{\L}(k^n))=\Frac({\mathcal O}_{\L}({k^*}^n)).$$
\end{enumerate}
\end{sdefin}

{\bf Remarque.}
Lorsque tous les coefficients de commutation $\l_{i,j}$ sont puissance d'un m\^eme $q\in k^*$ on parle de {\it cas uniparam\'etr\'e}.

\medskip

Nous d\'efinissons maintenant les alg\`ebres polynomiales mixtes crois\'ees, qui  nous serviront d'alg\`ebres de r\'ef\'erence pour le probl\`eme de Gelfand-Kirillov ``mixte'' (pour des exemples et des propri\'et\'es homologiques de ces alg\`ebres, voir  \cite{Rhhqw}).

\begin{sdefin} \label{vrdefqw}
Soient 
deux entiers $n\geq 1$ et $0\leq r\leq n$, et une matrice $\Lambda=(\lambda_{i,j})_{1\leq i,j\leq n}\in M_n(k^*)$ multiplicativement antisym\'etrique. 
On appelle alg\`ebre polynomiale mixte crois\'ee param\'etr\'ee par $n,r$ et $\L$, not\'ee $S_{n,r}^{\Lambda}(k)$, l'alg\`ebre engendr\'ee sur $k$ par $n+r$ g\'en\'erateurs $y_1,\ldots,y_n,x_1,\ldots,x_r$ avec les relations :
$$
\begin{array}{cll}
y_iy_j  = \l_{i,j}y_jy_i & \textrm{ si } & 1\leq i,j\leq n,\\
x_iy_i  =  y_ix_i+1 & \textrm{ si } & 1\leq i\leq r,\\
x_iy_j  =  \l_{i,j}^{-1}y_jx_i & \textrm{ si } & 1\leq i\leq r,\ 1\leq j\leq n,\ i\neq j,\\
x_ix_j  =  \l_{i,j}x_jx_i & \textrm{ si } & 1\leq i,j\leq r.\\
\end{array}
$$
\end{sdefin}

Remarquons que  la sous-alg\`ebre de $S_{n,r}^{\Lambda}(k)$ engendr\'ee par $y_1,\ldots,y_n$
est $k$-i\-so\-mor\-phe \`a l'espace quantique ${\mathcal O}_{\Lambda}(k^n)$,
et si $r\geq 1$, la sous-alg\`ebre de $S_{n,r}^{\Lambda}(k)$ engendr\'ee par
$x_i$ et $y_i$ est  $k$-isomorphe pour tout $1\leq i\leq r$ \`a l'alg\`ebre de Weyl $A_1(k)$.

Par ailleurs il est clair que $S_{n,r}^{\L}(k)$ est une extension it\'er\'ee de Ore, noeth\'erienne et int\`egre.
On peut donc consid\'erer son corps de fractions.

\begin{sdefin}
Soient 
deux entiers $n\geq 1$ et $0\leq r\leq n$, et une matrice $\Lambda\in M_n(k^*)$ multiplicativement antisym\'etrique. 
On appelle corps de fonctions rationnelles mixte crois\'e param\'etr\'e par $n,r$ et $\L$ le corps $\Frac(S_{n,r}^{\L}(k))$.
\end{sdefin}

La GK-dimension et le GK-degr\'e de transcendance sont deux invariants dimensionnels d\'efinis dans \cite{GK} afin de classer les alg\`ebres et les corps de Weyl.
La proposition suivante montre que pour les alg\`ebres $S_{n,r}^{\Lambda}(k)$  ils sont \'egaux au nombre de g\'en\'erateurs de l'alg\`ebre.

\begin{spropo} \label{gkdsnrl}
Soient $n\geq 1 $ et $0\leq r\leq n$ deux entiers, et $\L\in M_n(k^*)$ une matrice multiplicativement antisym\'etrique.
Alors
$$\gkdpr(S_{n,r}^{\L}(k))=\gktdpr(\Frac(S_{n,r}^{\L}(k)))=n+r.$$
\end{spropo}
\begin{demo}
Pour la dimension de Gelfand-Kirillov de l'extension it\'er\'ee de Ore $S_{n,r}^{\L}(k)$, on utilise r\'ecursivement le lemme 2.2 de \cite{HK}.
Le calcul du GK-degr\'e de transcendance du corps de fractions r\'esulte alors du th\'eor\`eme 7.3 de \cite{Z}.
Pour plus de d\'etails sur le calcul de ces invariants dans le contexte de cet article nous renvoyons au chapitre 1 de \cite{these}.
\end{demo}

\medskip

{\bf Notation.}
Pour $n=2$, la matrice $\L$ est param\'etr\'ee par un seul scalaire $\l=\l_{1,2}$.
On notera parfois $S_{2,r}^{\l}(k)$ pour $S_{2,r}^{\L}(k)$.

\medskip

{\bf Remarque.}
Comme on l'a d\'ej\`a not\'e dans \cite{Rhhqw}, les  alg\`ebres $\Alr$ peuvent \^etre vues comme des alg\`ebres d'op\'erateurs diff\'erentiels tordus sur un espace quantique param\'etr\'e par $\L$.

\subsection{Exemples et terminologie.} \label{13}

{\bf \ref{13}.1} {\sl Cas purement quantique.}
Dans le cas o\`u $r=0$, l'alg\`ebre $S_{n,0}^{\L}(k)$ n'est autre que  l'espace affine quantique ${\mathcal O}_{\L}(k^n)$.

\smallskip

{\bf \ref{13}.2} {\sl Cas purement classique.}
Dans le cas  o\`u tous les $\l_{i,j}$ valent 1, l'alg\`ebre $S_{n,r}^{(1)}(k)$ est l'alg\`ebre de Weyl  $A_r(k[y_{r+1},\ldots,y_n])=A_{r,t}(k)$ avec $t=n-r$. Lorsque  $r=n$ on retrouve l'alg\`ebre de Weyl usuelle $A_n(k)$.

\smallskip

{\bf \ref{13}.3} {\sl Cas semi-classique.}
Dans le cas $n=r$, sans hypoth\`ese sur la matrice $\L$,
les alg\`ebres $S_{n,n}^{\L}(k)$ sont des cas par\-ti\-cu\-li\-ers d'alg\`ebres de Weyl quantiques multiparam\'etr\'ees $A_n^{\qb,\L}(k)$ (voir \cite{AD1}, \cite{BG}, \cite{GZ}, \cite{Jslqwa},... ou la section 7 du pr\'esent article), correspondant au cas o\`u les param\`etres de quantification $\qb=(q_1,\ldots,q_n)$ valent tous 1. 
En particulier, l'alg\`ebre $S_{2,2}^{\L}(k)=A_2^{(1,1),\L}(k)$ est d\'ej\`a consid\'er\'ee dans \cite{AD3} (voir aussi le d\'ebut de la section 5 de cet article).
Les alg\`ebres $S^{\L}_{n,n}(k)$ sont simples, et de centre $k$ (voir \cite{GZ}, et le chapitre 5 de \cite{these}).
Il est d\'emontr\'e   dans \cite{Rhhqw} que  l'homologie et la cohomologie de Hochschild des $S^{\L}_{n,n}(k)$ sont identiques \`a celles de l'alg\`ebre de Weyl $A_n(k)$.


{\bf \ref{13}.4}
Beaucoup d'exemples d'alg\`ebres li\'ees aux groupes quantiques ont des localisations communes avec des alg\`ebres $\Alr$.
C'est bien s\^ur le cas de toutes les situations  o\`u l'analogue quantique de la conjecture de Gelfand-Kirillov est v\'erifi\'e (avec alors $r=0$).
C'est le cas de l'alg\`ebre \'etudi\'ee par D.A. Jordan dans \cite{Jed}, admettant pour corps de fractions le corps $\Frac(S_{2,1}^{q}(k))$, et  de certaines extensions it\'er\'ees de Ore (voir  le chapitre 1 de \cite{these}).

\smallskip

{\bf \ref{13}.5}
Les alg\`ebres polynomiales mixtes crois\'ees en dimension 2 sont l'alg\`ebre commutative $k[y_1,y_2]$ et le plan quantique $k_{\l}[y_1,y_2]$ pour $n=2$ et $r=0$, et l'alg\`ebre de Weyl $A_1(k)$ pour $n=r=1$.
On d\'eduit alors des r\'esultats de \cite{AVV} (voir aussi  la proposition 3.2 de \cite{AD1}) que toute extension de Ore it\'er\'ee en 2 variables est rationnellement \'equivalente \`a une alg\`ebre polynomiale mixte crois\'ee.
Il est facile de v\'erifier qu'en dimension 3 les diff\'erentes alg\`ebres polynomiales crois\'ees sont :
l'alg\`ebre commutative $k[x,y,z]$ et les espaces quantiques ${\mathcal O}_{\L}(k^3)$ pour $n=3$ et $r=0$, et  l'alg\`ebre de Weyl $A_{1,1}(k)=A_1(k[z])$ et les alg\`ebres  $S_{2,1}^{\L}(k)$ pour $n=2$ et $r=1$.

\subsection{Repr\'esentation graphique de certaines alg\`ebres.}

On reprend ici (en l'am\'enageant au cadre multiparam\'etr\'e) une convention introduite en \cite{AD3} pour illustrer les relations de commutation entre les g\'en\'erateurs de certains types d'extensions it\'er\'ees de Ore.

Soit $R$ une extension it\'er\'ee de Ore en $N$ g\'en\'erateurs $x_1,\ldots,x_N$.
On suppose que pour tous $1\leq i,j\leq N$ les relations entre g\'en\'erateurs sont de l'un des 2 types suivants:
$[x_i,x_j]=1$ ou $x_ix_j=\l_{i,j}x_jx_i$ avec $\l_{i,j}\in k^*$.
On repr\'esente alors $R$ par un graphe \`a $N$ sommets index\'es 
par $x_1,\ldots,x_N$ et tels que deux sommets quelconques $x_i$ et $x_j$ sont reli\'es par une ar\^ete orient\'ee ``color\'ee'' de l'une des deux fa\c cons suivantes:
\begin{itemize}
\item ar\^ete de Weyl (\'epaisse) :

\input{acc.pstex_t}\esp
signifie qu'on a la relation $[x_i,x_j]=1$;\\[.1cm]
\item ar\^ete quantique (mince, pond\'er\'ee) :

\input{aqq.pstex_t}\esp
 signifie qu'on a la relation $x_ix_j=\l_{i,j}x_jx_i$.
\end{itemize}

Pour all\'eger, on convient de ne pas faire figurer sur le graphe les ar\^etes quantiques de poids 1.
On a donc la convention suppl\'ementaire suivante :
\begin{itemize}
\item absence d'ar\^ete :

\input{padar.pstex_t}\esp
signifie qu'on a la relation $x_ix_j=x_jx_i$.
\end{itemize}

Ainsi l'alg\`ebre  de Weyl  ${A}_{n,t}(k)=A_n(k[z_1,\ldots,z_t])$ est repr\'esent\'ee par le graphe :
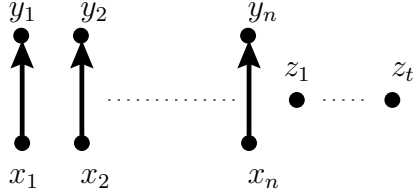
\begin{center}
\input{cw.pstex_t}
\end{center}

 Un espace quantique multiparam\'etr\'e ${\mathcal O}_{\L}(k^n)$ est repr\'esent\'e par un graphe du type :
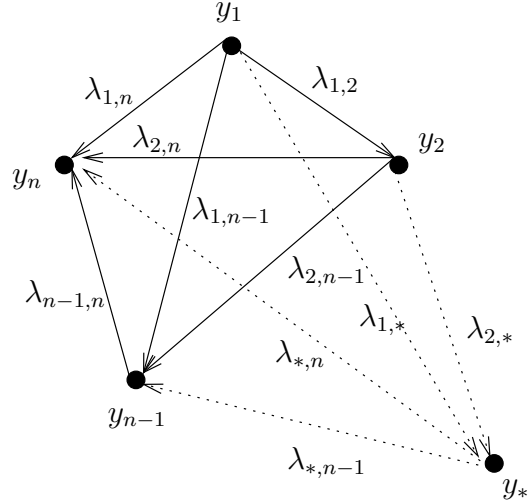
\begin{center}
\input{carre_q.pstex_t}
\end{center}

{\bf Remarque.}
A ce stade, de tels graphes ne sont que des illustrations permettant de pr\'esenter de fa\c con synth\'etique et parlante un certain nombre de relations de commutation.
La question de savoir reconna\^\i tre, parmis tous les graphes que l'on peut ainsi construire, ceux qui correspondent effectivement \`a des extensions it\'er\'ees de Ore en $N$ variables, ainsi que le probl\`eme de classifier les diff\'erents graphes repr\'esentant des alg\`ebres isomorphes sont trait\'es de fa\c con syst\'ematique \`a la section suivante.

\section{Graphes $qW$ et alg\`ebres polynomiales associ\'ees.} \label{secgraf}

Nous montrons dans cette section que les alg\`ebres polynomiales mixtes
crois\'ees sont la r\'eponse \`a la question (pos\'ee dans \cite{CGG}) de savoir
quelles classes d'alg\`ebres  \`a \'equivalence rationnelle pr\`es on obtient
\`a partir de pr\'esentations par g\'e\-n\'e\-ra\-teurs et relations ne faisant intervenir que des relations de type ``quantique'' ($XY=qYX$) ou ``de Weyl'' ($XY=YX+1$).

\subsection{Notion de graphe $qW$.}

\begin{sdefin}
On appelle graphe $qW$ la donn\'ee de $N$ sommets $X_1,\ldots$, $X_N$ deux \`a deux distincts reli\'es entre eux par des ar\^etes orient\'ees pond\'er\'ees des deux types suivants, arbitrairement (pour l'instant) d\'enomm\'ees :
\begin{itemize}
\item ar\^ete de Weyl de poids $p_{j,i} \in {\mathbb Z}$ :

\begin{center}
\input{ac.pstex_t}
\end{center}

\item ar\^ete quantique de poids $\l_{j,i} \in k^*$ : 

\begin{center}
\input{aq.pstex_t}
\end{center}

\end{itemize}
tels que :
\begin{enumerate}
\item aucune ar\^ete ne relie un sommet \`a lui-m\^eme;
\item
une ar\^ete au plus relie  deux sommets quelconques.
\end{enumerate}
\end{sdefin}

Lorsqu'une ar\^ete (quantique ou de Weyl) va du sommet $X_j$ vers le sommet $X_i$  on dit que $X_j$ est la \emph{source} de l'ar\^ete, et que $X_i$ est le \emph{but} de l'ar\^ete.

Soit $\Gamma$ un graphe $qW$ quelconque.
Il est clair que si l'on extrait de  $\Gamma$ un certain nombre de sommets et \emph{toutes} les ar\^etes qui relient ces sommets, on obtient encore un graphe $qW$, qu'on appelle un sous-graphe $qW$ de $\Gamma$.

\begin{sdefin} \label{defeqqw}
 On dit que deux graphes $qW$ sont \'equivalents si l'on passe de l'un \`a l'autre par un nombre fini d'op\'erations des types suivants, ou leurs inverses :
\begin{enumerate}
\item on place une ar\^ete de Weyl de poids 0 
\begin{center}
\input{ac0.pstex_t}
\end{center}
entre deux sommets reli\'es par aucune ar\^ete 
\begin{center}
\input{acomm.pstex_t}
\end{center}
\item on place une ar\^ete quantique de poids 1 
\begin{center}
\input{aq1.pstex_t} 
\end{center}
entre deux sommets reli\'es par aucune ar\^ete 
\begin{center}
\input{acomm.pstex_t}
\end{center}
\item on remplace une ar\^ete de Weyl de poids $p_{j,i}$ de $X_j$ vers $X_i$ 
\begin{center}
\input{ac.pstex_t} 
\end{center}
par une ar\^ete de Weyl de poids $-p_{j,i}$ de $X_i$ vers $X_j$ 
\begin{center}
\input{acinv.pstex_t}
\end{center}
\item on remplace une ar\^ete quantique de poids $\l_{j,i}$ de $X_j$ vers $X_i$ 
\begin{center}
\input{aq.pstex_t}
\end{center}
 par une ar\^ete quantique de poids $\l^{-1}_{j,i}$ de $X_i$ vers $X_j$ 
\begin{center}
\input{aqinv.pstex_t}
\end{center}
\end{enumerate} 
\end{sdefin}

\begin{sdefin} \label{defad}
 Un graphe $qW$ est dit {\rm admissible} si tous ses sous-graphes \`a trois sommets sont \'equi\-valents \`a l'un des graphes suivants :

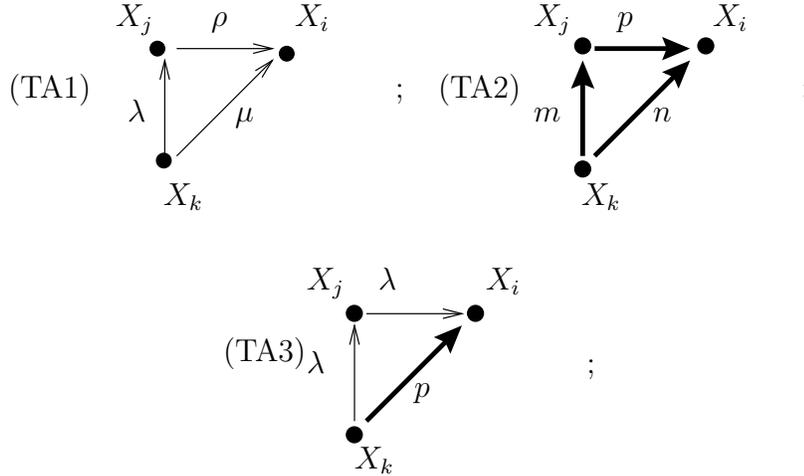
\begin{center}
\input{tradm1.pstex_t}
\end{center}

\begin{center}
\input{tradm2.pstex_t}
\end{center}

o\`u $m,n,p \in {\mathbb Z}$, et $\l,\mu,\rho \in k^*$.
\end{sdefin}

On parlera parfois de ``triangle admissible'' pour de tels sous-graphes.

\subsection{Alg\`ebre associ\'ee \`a un graphe $qW$.}

\begin{sdefin}
Soit $\G$ un graphe $qW$ \`a $N$ sommets $X_1,\ldots,X_N$. On appelle alg\`ebre associ\'ee \`a $\G$, et on note ${\mathcal A}(\G)$ le quotient de l'alg\`ebre  libre en $N$ g\'en\'erateurs $X_1,\ldots,X_N$ par l'id\'eal engendr\'e par les \'el\'ements $\{R_{j,i}\}_{(j,i)\in{\mathcal I}}$, o\`u ${\mathcal I}\subset \{1,\ldots,n\}^2$, 
 d\'efinis par :
\begin{itemize}
\item $R_{j,i} = X_jX_i - X_iX_j -p_{j,i}$  s'il existe une ar\^ete de Weyl de poids $p_{j,i}$ de  $X_j$ vers  $X_i$;
\item $R_{j,i}=X_jX_i -\l_{j,i}X_iX_j$ s'il existe une ar\^ete quantique de poids $\l_{j,i}$ de $X_j$ vers $X_i$;
\item $R_{j,i}=X_jX_i - X_iX_j$ si $X_j$ et $X_i$ ne sont reli\'es par aucune ar\^ete.
\end{itemize}
\end{sdefin}

Dans la suite, $X_i$ d\'esignera  indiff\'eremment un sommet d'un graphe $\G$ ou un g\'en\'erateur de l'alg\`ebre $\AG$.
Il est clair que deux graphes $qW$ \'equivalents ont des alg\`ebres associ\'ees $k$-isomorphes.
Par ailleurs, une alg\`ebre $S_{n,r}^{\L}(k)$ appara\^\i t clairement comme une alg\`ebre associ\'ee \`a un graphes $qW$, polynomiale en ces g\'en\'erateurs.
On d\'emontre dans la suite de cette section qu'\`a isomorphisme pr\`es les alg\`ebres polynomiales mixtes crois\'ees sont les seules telles alg\`ebres.

\begin{slemme}  \label{defadmqw}
Soit $\G$ un graphe $qW$ \`a $N$ sommets $X_1,\ldots,X_N$.
On suppose que dans l'alg\`ebre $\AG$, les g\'en\'erateurs $X_1,\ldots,X_N$ sont $k$-lin\'eairement ind\'ependants.
 Soient $X_i$ et $X_j$ deux g\'en\'erateurs de ${\mathcal A}(\G)$ v\'erifiant $X_iX_j = X_jX_i +p$, avec $p \in {\mathbb Z} \setminus \{0\}$. Soit $X_k$ un autre g\'en\'erateur. On a les deux propri\'et\'es suivantes :
\begin{enumerate}
\item si $X_iX_k-X_kX_i=m \in {\mathbb Z}$, alors $X_jX_k - X_kX_j =n \in {\mathbb Z}$ ;
\item si $X_iX_k= \l X_kX_i$, avec $\l \in k\setminus \{0,1\}$, alors $X_jX_k=\l^{-1} X_kX_j$.
\end{enumerate}
\end{slemme}
\begin{demo}
{1}. Supposons que $X_kX_j = \l X_jX_k$. Puisque $\AG$ est associative, on a : $(X_kX_j)X_i=X_k(X_jX_i)$. Or 
$$(X_kX_j)X_i=\l X_jX_kX_i=\l X_j(X_iX_k -m)=\l X_iX_jX_k -\l pX_k -\l mX_j.$$
Par ailleurs, 
$$X_k(X_jX_i)= X_k(X_iX_j -p)=X_iX_kX_j - mX_j -pX_k= \l X_iX_jX_k -mX_j-pX_k.$$
 Par hypoth\`ese sur ${\mathcal A}(\G)$, les g\'en\'erateurs $X_j$ et $X_k$ sont ind\'ependants,  on doit donc avoir $\l p=p$ et $\l m=m$. Puisque $p \neq 0$, on conclut $\l =1$.


{2}. En inversant les r\^oles de $X_i$ et $X_j$,  l'hypoth\`ese 
$X_iX_j-X_jX_i=p\in{\mathbb Z}\setminus\{0\}$ et  le point {1} impliquent 
 qu'on ne peut pas avoir \`a la fois $X_iX_k= \l X_kX_i$, avec $\l \in k\setminus \{0,1\}$, et $X_jX_k - X_kX_j =n \in {\mathbb Z} \setminus \{0\}$.
On a donc : 
$X_jX_k = \mu X_kX_j$, avec $\mu\in k^*$. Alors 
$$(X_kX_j)X_i=\mu^{-1} X_jX_kX_i=\mu^{-1}\l^{-1} X_jX_iX_k =\mu^{-1}\l^{-1} X_iX_jX_k - \mu^{-1}\l^{-1} pX_k.$$
Par ailleurs, 
$$X_k(X_jX_i)= X_k(X_iX_j -p)=\l^{-1}X_iX_kX_j  -pX_k= \l^{-1} \mu^{-1}X_iX_jX_k -pX_k.$$
 A nouveau de $(X_kX_j)X_i=X_k(X_jX_i)$ on d\'eduit que  $ \mu^{-1}\l^{-1} p=p$. Puisque $p \neq 0$, on conclut $\l \mu=1$.
\end{demo}

\begin{spropo} \label{propadm}
Soit $\G$ un graphe $qW$ de sommets $X_1,\ldots,X_N$. Alors les deux assertions suivantes sont \'equivalentes :
\begin{enumerate}
\item les mon\^omes $\{X_1^{\a_1}\ldots X_N^{\a_N}\}_{\a=(\a_1,\ldots,\a_N) \in {\mathbb N}^N}$ forment une base de $k$-espace vectoriel de l'alg\`ebre ${\mathcal A}(\G)$;
\item le graphe $\G$ est admissible.
\end{enumerate}
\end{spropo}

\begin{demo}
A \'equivalence pr\`es, les sous-graphes \`a trois sommets d'un graphe $qW$  quelconque $\G$ sont de l'un des types suivants :

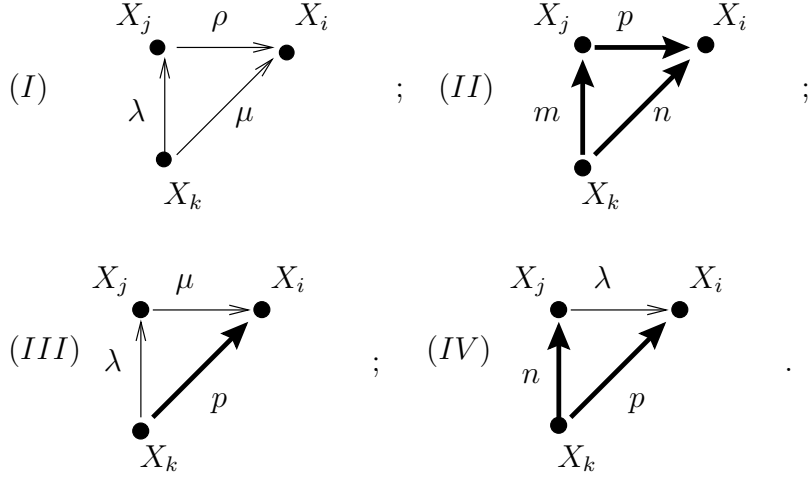
\begin{center}
\input{trposs.pstex_t}
\end{center}

On remarque que les triangles (I) et (II) sont toujours admissibles.
Supposons que le point (1) soit vrai. 
Alors le lemme \ref{defadmqw} s'applique, et permet de montrer que 
le cas (III) ne peut avoir lieu que si $p=0$ (il se ram\`ene alors au cas (I)) ou si $\l\mu=1$ (ce qui donne un triangle admissible de type (TA3)).
Le cas (IV) ne peut avoir lieu que si $\l=1$ (il se ram\`ene alors au cas (II)) ou si $n=p=0$ (il se ram\`ene alors au cas (I)).

R\'eciproquement,  supposons que $\G$ soit admissible. On montre que les mon\^omes $\{X_1^{a_1}\ldots X_N^{a_N}\}$ forment une base de $k$-espace vectoriel de ${\mathcal A}(\G)$ gr\^ace au lemme du diamant. 
Nous renvoyons le lecteur au chapitre I.11 de \cite{BG} pour une
pr\'e\-sen\-ta\-ti\-on de ce r\'esultat
adapt\'ee au cadre de cet article, et pour les d\'efinitions des notions de syst\`eme de r\'eduction, d'ambigu\"\i t\'e d'inclusion et d'ambigu\"\i t\'e de superposition.
A \'equivalence pr\`es, on peut r\'eorienter  les fl\`eches de $\G$ de telle sorte qu'elles soient toutes orient\'ees de $X_j$ vers $X_i$ pour $j>i$.
En suivant les notations du chapitre I.11 de \cite{BG}, l'ordre longueur-lexicographique (c'est-\`a-dire qu'on ordonne les mots d'abord par leur longueur, puis on ordonne les mots de longueur \'egale par l'ordre lexicographique usuel)  sur les mon\^omes est compatible avec le syst\`eme de r\'eduction 
$$S=\{(X_jX_i,X_jX_i-R_{j,i}),\ 1\leq i<j\leq N\},$$
et il n'y a pas d'ambigu\"\i t\'e d'inclusion, car les $X_jX_i$ sont tous des mots distincts de la m\^eme longueur.
On peut alors avoir 8 types d'ambigu\"\i t\'es de superposition, suivant les relations que v\'erifient $X_i,X_j,X_k$ avec $i<j<k$.
La v\'erification de la r\'esolubilit\'e de toutes ces ambigu\"\i t\'es de superposition est  facile sous l'hypoth\`ese que $\G$ est admissible.
A titre d'exemple,
soient  $i<j<k$ tels que 
$$R_{j,i}=X_jX_i-X_iX_j-1,\ R_{k,i}=X_kX_i-\l X_iX_k,\ R_{k,j}=X_kX_j-\l^{-1}X_jX_k.$$
Alors par r\'esolutions successives 

\begin{center}
$(X_kX_j)X_i$ donne $\l^{-1}X_jX_kX_i$, qui donne $\l^{-1}X_j\l X_iX_k$, qui donne $(X_iX_j+1)X_k$.
\end{center}

D'autre part, 

\begin{center}
$X_k(X_jX_i)$ donne $X_k(X_iX_j+1)$, qui  donne $\l X_iX_kX_j+X_k$, qui  donne $\l X_i\l^{-1}X_jX_k+X_k$,
\end{center}

 et on a bien $(X_iX_j+1)X_k=X_iX_jX_k+X_k$.
Les autres ambigu\"\i t\'es de superposition se traitent  de m\^eme,  et sont laiss\'ees au lecteur.
\end{demo}

\begin{scorol}
Soit $\G$ un graphe $qW$ admissible avec $N$ sommets $X_1,\ldots$, $X_N$. Alors :
 \begin{enumerate}
\item
$\AG$ est une extension it\'er\'ee de Ore $k[X_1][X_2;\sigma_2,\delta_2]\ldots [X_N;\sigma_N,\delta_N]$,
o\`u les automorphismes et $\sigma$-d\'erivations traduisent les relations entre les g\'en\'erateurs $X_i$ cod\'ees par le graphe $\G$;
\item
$\gkdpr(\AG)=\gktdpr(\Frac(\AG))=N$.
\end{enumerate}
\end{scorol}
\begin{demo}
1. Le fait que $\AG$ soit une extension it\'er\'ee de Ore  d\'ecoule de l'existence d'une base de $k$-espace vectoriel montr\'ee dans la proposition pr\'ec\'edente, et des relations de commutation entre les g\'en\'erateurs, qu'on traduit \`a l'aide des automorphismes $\sigma_i$ et des $\sigma_i$-d\'erivations $\delta_i$. \\[.2cm]
2. Comme \`a la proposition \ref{gkdsnrl}, on utilise r\'ecursivement le lemme 2.2 de \cite{HK}, puis le th\'eor\`eme 7.3 de \cite{Z}.
\end{demo}

\medskip

{\bf Remarque.} Les alg\`ebres polynomiales crois\'ees $S_{n,r}^{\L}(k)$ de la section \ref{defsnrl} sont, de par leur d\'efinition m\^eme, des exemples d'alg\`ebres $\AG$.
Le but de la fin de cette section  est de d\'emontrer que r\'eciproquement, toute alg\`ebre $\AG$ associ\'ee \`a un graphe admissible est isomorphe \`a une alg\`ebre $S_{n,r}^{\Lambda}(k)$.

\subsection{R\'eduction des graphes $qW$.} \label{redgrsec}

Plusieurs graphes $qW$ peuvent \'evidemment \^etre associ\'es \`a des alg\`ebres isomorphes. Ainsi par exemple 
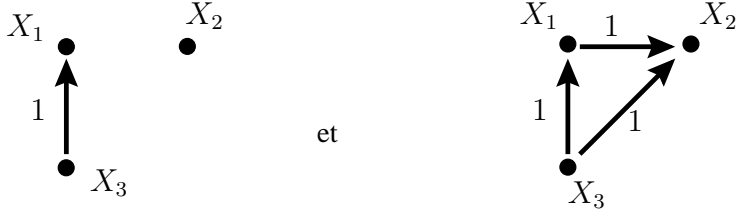
\begin{center}
\input{cw11.pstex_t}
\end{center}

  repr\'esentent tous les deux l'alg\`ebre  de Weyl $A_{1,1}(k)$ (dans la seconde alg\`ebre il suffit de remplacer $X_2$ par un autre g\'en\'erateur  : $X_1-X_2-X_3$ pour retrouver les relations usuelles de $A_{1,1}(k)$).
On cherche donc \`a r\'eduire les graphes $qW$ \`a une forme canonique, qui correspondra exactement aux alg\`ebres polynomiales mixtes crois\'ees.

\begin{sdefin} \label{gengr}
Soit $\G$ un graphe $qW$ admissible \`a $N$ sommets. On appelle  famille de {\rm g\'en\'erateurs  graphiques} de l'alg\`ebre $\AG$ toute  famille $(Y_1,\ldots,Y_N)$ d'\'el\'ements de $\AG$, v\'erifiant les points suivants :
\begin{enumerate}
\item les mon\^omes $\{Y_1^{\a_1}\ldots Y_N^{\a_N}\}_{\a\in{\mathbb N}^N}$ forment une base de $k$-espace vectoriel de $\AG$;
\item les $Y_i$ v\'erifient seulement des relations de type quantique ($Y_i Y_j = \l_{i,j} Y_jY_i$, $\l_{i,j}\in k^*$) ou de type Weyl pond\'er\'e ($Y_iY_j=Y_jY_i +p_{i,j}$, $p_{i,j}\in{\mathbb Z}$).
\end{enumerate}
\end{sdefin}

Par d\'efinition, les sommets d'un graphe $qW$ admissible $\Gamma$ sont des g\'en\'erateurs graphiques de $\AG$.
Et si $(Y_1,\ldots,Y_N)$ est une famille de g\'en\'erateurs graphiques d'une alg\`ebre $\AG$ associ\'ee \`a  graphe $qW$ admissible $\Gamma$,
on peut  construire un graphe $qW$ admissible $\G'$,  de sommets $Y_1,\ldots,Y_N$ et dont les ar\^etes codent les relations entre les $Y_i$ dans $\AG$. Evidemment on a alors $\AG = {\mathcal A}(\G')$.
On dira dans la suite que $\G'$ est le graphe $qW$ admissible \emph{associ\'e} \`a $(Y_1,\ldots,Y_N)$.

\medskip

{\bf Notations.}
Pour un graphe $qW$ fix\'e $\G$ quelconque, on note ${\mathcal W}(\G)$ l'ensemble des sommets de $\G$ qui sont sommet
 d'au moins une ar\^ete de Weyl de poids non nul de $\G$.
 On note ${V}(\G)$ l'ensemble des paires $\{X_i,X_j\}$ de sommets distincts de $\G$, 
tels que $X_i$ et $X_j$ sont reli\'es par une ar\^ete de Weyl de poids non nul (dans un sens ou dans l'autre), et tels que ni $X_i$ ni $X_j$ ne sont sommets d'une autre ar\^ete de Weyl de poids non nul.
 Enfin, on note ${\mathcal V}(\G)$ l'ensemble des sommets de $\G$ \'el\'ements d'une paire $\{X_i,X_j\}$ appartenant \`a  $ V(\G)$. On remarque que ${\mathcal V}(\G) \subset {\mathcal W}(\G)$.

\medskip

Le but de ce qui suit est de montrer que, {\it via} un changement de g\'en\'erateurs graphiques, on peut toujours se ramener \`a ${\mathcal V}(\G) ={\mathcal W}(\G)$.
Cette r\'eduction est obtenue au lemme \ref{redgr} gr\^ace \`a la suite de lemmes suivants.

\begin{slemme}  \label{l1}
Soit $\G$ un graphe $qW$ admissible. Soient $X,Y,Z$ trois sommets de $\G$, tels que $XY=YX+n$ et $XZ=ZX+p$, avec $n,p\in{\mathbb Z}\setminus\{0\}$. Notons $T_1,\ldots,T_k$ les autres sommets de $\G$. 
Dans l'alg\`ebre $\AG$, consid\'erons $Y' = \a Y+ \b Z$, o\`u
 $\a,\b \in {\mathbb Z}$, $\a\neq 0$. Alors $X,Y',Z,T_1, \ldots,T_k$ sont des g\'en\'erateurs graphiques de $\AG$. Plus pr\'ecis\'ement :
\begin{itemize}
\item si $YT_i=\l T_iY$, avec $\l \neq 1$, alors $Y'T_i=\l T_iY'$.
\item si $YT_i=T_iY +m$, avec $m\in{\mathbb Z}$, alors il existe $r\in{\mathbb Z}$ tel que $ZT_i=T_iZ + r$, et $Y'T_i=T_iY' +\a m +\b r$.
\end{itemize}
\end{slemme}
\begin{demo}
La famille $(X,Y',Z,T_1, \ldots,T_k)$ v\'erifie le point 1  de la d\'efinition  \ref{gengr}. 
Par ailleurs si $YT_i=\l T_iY$, avec $\l \neq 1$, alors $XT_i=\l^{-1}T_iX$ (lemme \ref{defadmqw}), et $ZT_i=\l T_iZ$, d'o\`u  $Y'T_i=\l T_iY'$. De m\^eme  si $YT_i=T_iY +m$,  on ne peut pas avoir $ZT_i=\l T_iZ$ pour $\l \neq 1$, et le reste du lemme se v\'erifie facilement.
\end{demo}

\begin{slemme} \label{chouchy}
Soit $\G$ un graphe $qW$ admissible de sommets $X_1,\ldots,X_N$.
Si  ${\mathcal V}(\G)\neq {\mathcal W}(\G)$, alors il existe une famille $(Y_1,\ldots,Y_N)$  de g\'en\'erateurs graphiques de $\AG$,  tels que le graphe   $qW$ admissible $\G'$ associ\'e \`a $(Y_1,\ldots,Y_N)$ v\'erifie :
$$ \sharp {\mathcal V}(\G') > \sharp {\mathcal V}(\G), \ \textrm{ et } \ \sharp {\mathcal W}(\G') \leq \sharp {\mathcal W}(\G).$$
\end{slemme}
\begin{demo}
La d\'emonstration se fait
 en plusieurs \'etapes. Soit  $X \in {\mathcal W}(\G) \setminus {\mathcal V}(\G)$.
Notons $\G_o$ le sous-graphe de $\G$ constitu\'e des sommets de $\G$ qui sont dans ${\mathcal V}(\G)$.

{\it Etape 1}.
Supposons $X$  reli\'e par $r$ ar\^etes de Weyl de poids $p_1,\ldots, p_r$ non nuls respectivement \`a des  sommets $S_1,\ldots, S_r$ de $\G$. Notons $Z_1,\ldots, Z_t$ les sommets restants de $\G$. 
Parce que $X\not\in {\mathcal V}(\G)$, les sommets $S_1,\ldots,S_r$ ne sont pas non plus dans ${\mathcal V}(\G)$, et  ${\mathcal V}(\G)$ est inclus dans  $\{Z_1,\ldots,Z_t\}$.
 Si $r=1$ on pose $\wt\G = \G$ et on passe directement \`a la deuxi\`eme \'etape. Sinon, on a partiellement le sch\'ema suivant :

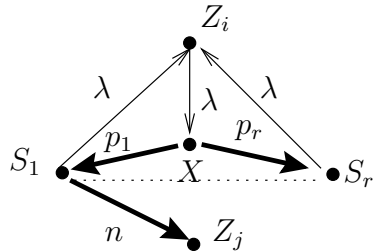
\begin{center}
\input{etap1.pstex_t}
\end{center}

 On va montrer qu'\`a un changement de g\'en\'erateurs pr\`es, on peut supposer que $X$ n'est   sommet que d'une ar\^ete de Weyl de poids non nul. Pour cela on d\'emontre qu'il
 existe des \'el\'ements $S'_1,\ldots, S'_r$ de $A$, v\'erifiant les points suivants :
\begin{enumerate}
\item les \'el\'ements  $X, S'_1,\ldots, S'_r, Z_1,\ldots, Z_t$ sont des g\'en\'erateurs graphiques de $A$ ;
\item pour tout $i \in \{1,r\}$, on a $[X,S'_i]= p'_i$, et $p'_1$ divise $p'_i$ pour tout $i\geq 2$ tel que  $p'_i \neq 0$ ;
\item si on note $\what\G$ le graphe $qW$ admissible associ\'e \`a $(X,S'_1,\ldots,S'_r,Z_1,\ldots,Z_t)$, on a pour tout $j$ les implications : 
$$[Z_j\in {\mathcal W}(\what\G) \oif Z_j\in {\mathcal W}(\G)]\ \textrm{ et }\  [Z_j\in {\mathcal V}(\G) \iff Z_j\in {\mathcal V}(\what\G)].$$
\end{enumerate}
 Supposons d'abord que $r=2$. Si $p_1$ divise $p_2$ ou si $p_2$ divise $p_1$, quitte \`a r\'eindexer $S_1$ et $S_2$ il n'y a rien \`a faire.
Sinon soit $d={\rm pgcd} (p_1,p_2)$, de sorte que $d=up_1+vp_2$, avec $u$ et $v$ entiers non nuls.
On pose alors $S'_1=uS_1+vS_2$ et  $S'_2=S_2$.
 Le lemme \ref{l1} permet de conclure que les points (1) et  (2) ci-dessus sont bien v\'erifi\'es, et que  pour tout $j\leq t$, s'il existe $\l\neq 1$ tel que $S_1Z_j = \l Z_jX_1$, alors $S_2Z_j = \l Z_jS_2$.
 Il d\'ecoule de ceci  qu'en rempla\c cant $S_1$ par $S'_1$ on ne va ``cr\'eer'' de nouvelle ar\^ete  de Weyl entre $Z_j$ et $S'_1$ que si $Z_j$ \'etait d\'ej\`a reli\'e par une ar\^ete de Weyl \`a
  $S_1$ ou $S_2$.
 Ainsi on ne fait pas appara\^\i tre dans ${\mathcal W}(\what\G)$ des sommets $Z_j$ qui ne se trouvaient pas dans ${\mathcal W}(\G)$, et une paire de sommets $\{Z_i,Z_j\}$ dans $V(\G)$ reste dans $V(\what\G)$.
De l\`a d\'ecoule le point {3}. 

Supposons maintenant $r\geq 3$. On effectue l'op\'eration pr\'ec\'edente avec $S_{r-1}$ et $S_r$.
On obtient ainsi deux nouveaux g\'en\'erateurs $S_{r-1}'$ et $S_r'$ v\'erifiant les trois points ci-dessus. On recommence alors avec $S_{r-2}$ et $S_{r-1}'$, et ainsi de suite jusqu'\`a obtenir une famille $S'_1,\ldots,S'_r$ v\'erifiant   les conditions voulues.

\medskip

Le graphe $\what\G$ v\'erifiant (1), (2) et (3) \'etant ainsi construit, on d\'efinit  de nouveaux g\'en\'erateurs 
 $S''_k = S'_k - (p'_k / p'_1) S'_1$ dans $A$, pour tout $k\geq 2$, en notant que $p'_k/p'_1\in{\mathbb Z}$.
 Par construction  $X$ commute \`a $S''_2,\ldots,S''_r$. Comme ci-dessus, on a une nouvelle famille de g\'en\'erateurs graphiques, associ\'es \`a un nouveau graphe $\wt\G$. 
On note que, \`a cause du point (3) v\'erifi\'e par les $S'_i$,  on a $\sharp {\mathcal W}(\wt\G) \leq \sharp {\mathcal W}(\G)$ et ${\mathcal V}(\wt\G)$ contient $\G_o$. 

\medskip

{\it Etape 2}.
 On suppose donc que $X$ n'est sommet dans $\wt\G$ que d'une ar\^ete de Weyl, de poids non nul $p$. Soit $Y$ l'autre sommet de cette ar\^ete. 
Si, du fait des manipulations effectu\'ees \`a l'\'etape 1, $Y$ n'est plus rattach\'e \`a aucune autre ar\^ete de Weyl de poids non nul (cette situation ne pouvait pas \^etre la situation initiale puisqu'on a suppos\'e $X \not\in {\mathcal V}(\G)$ ), on pose $\G'=\wt\G$ et on passe \`a l'\'etape 3.

Sinon,  on a partiellement le sch\'ema suivant :

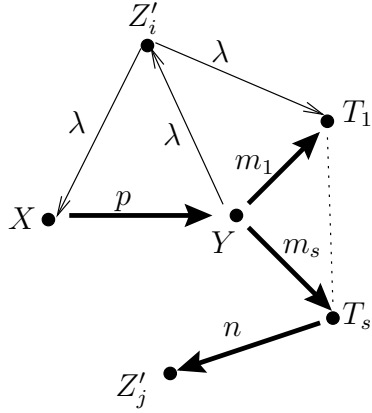
\begin{center}
\input{etap2.pstex_t}
\end{center}

  Quitte \`a remplacer $X$ par $p^{-1}X$, on peut supposer que $p=1$. Changeons la notation des sommets de $\wt\G$. Le sommet $Y$ est li\'e par des ar\^etes de Weyl de poids non nul \`a $X$ et \`a des sommets $T_1,\ldots,T_s$.
Quitte \`a remplacer $\wt\G$ par un graphe \'equivalent, on peut supposer que les ar\^etes de Weyl sont orient\'ees de $Y$ vers $T_1,\ldots,T_s$.
Notons alors $m_1,\ldots,m_s$ le poids de ces ar\^etes.
 Notons  $Z'_1,\ldots,Z'_r$ les sommets restants. Comme pr\'ec\'edemment ${\mathcal V}(\wt\G)$, et donc $\G_o$ sont inclus dans $\{Z_1,\ldots,Z_t\}$.
 On r\'eit\`ere alors  avec $Y$ ce qui a \'et\'e fait pour $X$ \`a la fin de l'\'etape 1, en posant $T'_i=T_i + m_iX$ pour tout $i$, ce qui donne une nouvelle famille $(X,Y,T'_1,\ldots,T'_s,Z_1,\ldots,Z_r)$ de g\'en\'erateurs graphiques pour $\AG$, tels que $X$ ne soit li\'e par une ar\^ete de Weyl qu'\`a $Y$ et $Y$ qu'\`a $X$. Notons $\G'$ le graphe $qW$ admissible qui leur est associ\'e.
Comme pr\'ec\'edemment, on a 
$$\sharp {\mathcal W}(\G') \leq \sharp {\mathcal W}(\wt\G),\ \textrm{ et }\ [Z'_j\in {\mathcal V}(\G') \iff Z'_j\in {\mathcal V}(\wt\G)].$$

\medskip

{\it Etape 3}. 
De tout ce qui pr\'ec\`ede il d\'ecoule que $\{X,Y\} \in V(\G')$, et  les sommets de $\G_o$ (c'est-\`a-dire les sommets initialement dans ${\mathcal V}(\G)$ ) sont encore dans ${\mathcal V}(\G')$. On a donc  $ \sharp {\mathcal V}(\G') \geq \sharp {\mathcal V}(\G) +2 > \sharp {\mathcal V}(\G)$.
 Par ailleurs on s'est assur\'e que l'on avait $\sharp {\mathcal W}(\G') \leq \sharp {\mathcal W}(\wt\G) \leq \sharp {\mathcal W}(\G)$. Ainsi le lemme est d\'emontr\'e.
\end{demo}

\begin{slemme}     \label{redgr}
Soit $\G$ un graphe admissible de sommets $X_1,\ldots,X_N$. Alors il existe une famille $(Y_1,\ldots,Y_N)$ de g\'en\'erateurs graphiques de ${\mathcal A} (\G)$,  tels que le graphe $qW$ admissible $\G'$ associ\'e \`a $(Y_1,\ldots,Y_N)$  v\'erifie ${\mathcal V}(\G') = {\mathcal W}(\G')$. Ainsi, tout sommet de $\G'$ est sommet d'{\rm au plus} une ar\^ete de Weyl de poids non nul.
\end{slemme}
\begin{demo}
On raisonne par r\'ecurrence sur le cardinal de ${\mathcal W}(\G) \setminus {\mathcal V}(\G)$, \`a l'aide du lemme \ref{chouchy}.
\end{demo}

\begin{stheoreme} \label{redgrpropo}
Soit $\G$ un graphe $qW$ admissible \`a $N$ sommets.
Alors il existe deux entiers $n\geq 1$ et $0\leq r\leq n$ v\'erifiant $N=n+r$, et une matrice multiplicativement antisym\'etrique $\L\in M_n(k^*)$ tels que $N=n+r$ et l'alg\`ebre  $\AG$  est $k$-isomorphe \`a l'alg\`ebre $\Alr$.
\end{stheoreme}
\begin{demo}
On peut sans restriction supposer que le graphe $\G$ v\'erifie les conclusions du lemme \ref{redgr}.
Il est clair par ailleurs (voir le point 1 de la d\'efinition \ref{defeqqw}) que l'on peut aussi supposer que $\G$  ne contient 
aucune ar\^ete de Weyl de poids nul.
Soit alors $r$ le nombre d'ar\^etes de Weyl de $\G$.
 Pour chaque $i\leq r$ notons $X_i$ et $Y_i$ la source et le but de la $i^{\smmath{\`eme}}$ ar\^ete de Weyl, et $p_{i}$ le poids de cette ar\^ete.
Par les conclusions du  lemme \ref{redgr} les sommets $X_1,Y_1,\ldots,X_r,Y_r$ sont tous distincts.
En rempla\c cant $X_i$ par $p_{i}^{-1} X_i$, on  peut supposer sans perte de g\'en\'eralit\'e que toutes les ar\^etes de Weyl sont de poids 1.
Posons $n=N-r$, et notons   $Y_{r+1}, \ldots,Y_n$ les sommets de $\G$ qui ne sont ni des $X_i$, ni des $Y_i$ pour $1\leq i\leq r$. 
L'alg\`ebre $\AG$ appara\^\i t alors comme l'alg\`ebre engendr\'ee sur $k$ par $X_1,Y_1,\ldots,X_r,Y_r,Y_{r+1}, \ldots,Y_n$
avec les relations donn\'ees en \ref{vrdefqw},  o\`u les coefficients $\l_{i,j}$  sont donn\'es par le poids de l'ar\^ete quantique liant $Y_i$ \`a $Y_j$ pour tous $1\leq i,j\leq n$.
On conclut que $\AG$ est $k$-isomorphe \`a $S_{n,r}^{\L}(k)$.
\end{demo}

\medskip

{\bf Remarque.}
Soit $\G$ un graphe $qW$ admissible \`a $N$ sommets $X_1,\ldots,X_N$. Notons $P(\G)=(p_{i,j})$ la matrice antisym\'etrique de $M_N({\mathbb Z})$ dont le $(i,j)^{\smmath{\`eme}}$ coefficient vaut 0  si $X_i$ et $X_j$ sont reli\'es par une ar\^ete quantique ou par aucune ar\^ete,  et le poids de l'ar\^ete de Weyl de $X_i$ vers $X_j$ (on peut toujours choisir l'orientation d'une ar\^ete de Weyl gr\^ace au point 3 de la d\'efinition \ref{defeqqw}) sinon.
On peut observer au cours des  d\'emonstrations des lemmes \ref{l1}, \ref{chouchy} et \ref{redgr},
 que les changements de g\'en\'erateurs graphiques effectu\'es \`a chaque \'etape correspondent \`a 
 des op\'erations \'el\'ementaires sur la matrice $P(\G)$.
On ne change donc pas son rang.
Puisque $P(\G')$ est clairement une matrice antisym\'etrique  de rang $2r$,
on retiendra  que l'entier $r$ intervenant dans le th\'eor\`eme \ref{redgrpropo} n'est autre que la moiti\'e du rang de la matrice $P(\G)$ associ\'ee au graphe de d\'epart (et $n$ est alors donn\'e par $n=N-r$, o\`u $N$ est le nombre de sommets de $\G$).

\medskip

Les entiers $n$ et $r$ du th\'eor\`eme \ref{redgrpropo} sont donc d\'etermin\'es de fa\c con unique et s'interpr\`etent non seulement en terme d'isomorphisme des alg\`ebres $S_{n,r}^{\L}(k)$, mais m\^eme d'\'equivalence rationnelle, comme on va le voir \`a la section \ref{secwdeg}.

\section{Premiers invariants des corps de fonctions rationnelles  mixtes crois\'es.}

\subsection{Pr\'esentation  eul\'erienne des corps $\Frac(S_{n,r}^{\Lambda}(k))$.}  \label{arres}

\subsubsection{Notations.}

Soient $n\geq 1$ et $0\leq r\leq n$ deux entiers, et $\L\in M_n(k^*)$ multiplicativement antisym\'etrique.
Dans l'alg\`ebre polynomiale mixte crois\'ee
 $\Alr$  d\'efinie en \ref{vrdefqw}, on d\'efinit les \'el\'ements  :
$$w_i=y_ix_i \ \textrm{ pour tout }\ 1\leq i\leq r.$$
On a alors  :
$$\begin{array}{l}
{[w_i,w_j]}=0\quad \forall  i,j\in \{1,\ldots,r\},\\
{[w_i,y_k]} = \delta_{i,k}y_k\quad  \forall  k \in\{1,\ldots, n\},\ \forall  i\in \{1,\ldots,r\}.
\end{array}$$
Notons $T_{n,r}^{\L}(k)$ la sous-alg\`ebre de $\Alr$ engendr\'ee par $w_1,\ldots,w_r,y_1,\ldots,y_r$, $y_{r+1}$,$\ldots$,$y_n$, et $\what T_{n,r}^{\L}(k)$ le localis\'e de $T_{n,r}^{\L}(k)$ en la partie multiplicative engendr\'ee par les \'el\'ements normalisants $y_1,\ldots,y_n$.
L'alg\`ebre $T_{n,r}^{\L}(k)$ [resp. $\what T_{n,r}^{\L}(k)$] s'interpr\`ete ais\'ement comme la sous-alg\`ebre de $\End({\mathcal O}_{\L}(k^n))$ [resp. de $\End({\mathcal O}_{\L}({k^*}^n))$] engendr\'ee par les op\'erateurs de multiplication \`a gauche par les g\'en\'erateurs $y_1,\ldots,y_n$ [resp. et leurs inverses] et les d\'eriv\'ees eul\'eriennes $d_1,\ldots,d_r$ (qui sont des vraies d\'erivations), d\'efinies par $d_i(y_k)=\d_{i,k}y_k$ pour tous $1\leq k\leq n$, $1\leq i\leq r$.

\medskip

Il est facile de v\'erifier que $T_{n,r}^{\L}(k)$ est l'extension de Ore it\'er\'ee  :
$$T_{n,r}^{\L}(k)=k[w_1,\ldots,w_r][y_1;\sigma_1]\ldots[y_n;\sigma_n],$$
o\`u chaque automorphisme $\sigma_i$ de $k[w_1,\ldots,w_r][y_1;\sigma_1]\ldots[y_{i-1};\sigma_{i-1}]$ est d\'efini par  :
$$\begin{array}{l}
\sigma_i(w_j)=w_j-\d_{i,j} \ \textrm{ pour tout }  1\leq j\leq r,\\
\sigma_i(y_j)=\l_{i,j}y_j \ \textrm{ pour tout }  1\leq j<i.
\end{array}$$
En particulier on a  :
$$T_{n,r}^{\L}(k)\subset S_{n,r}^{\L}(k)\subset \what T_{n,r}^{\L}(k),$$
et les alg\`ebres $T_{n,r}^{\L}(k)$, $\what T_{n,r}^{\L}(k)$  et $S_{n,r}^{\L}(k)$ sont rationnellement \'equivalentes  :
$$\begin{array}{c}
\Frac(T_{n,r}^{\L}(k))=\Frac(\what T_{n,r}^{\L}(k))=\Frac(S_{n,r}^{\L}(k))=\\
\noalign{\smallskip}
k(w_1,\ldots,w_r)(y_1;\sigma_1)\ldots(y_n;\sigma_n),\end{array}$$
o\`u $\sigma_1,\ldots,\sigma_n$ d\'esignent encore les prolongements des automorphismes ci-dessus au corps de fractions des alg\`ebres consid\'er\'ees.
Les corps de fonctions rationnelles mixtes crois\'es apparaissent ainsi comme des exemples des classes de corps introduits par V.A. Artamonov dans \cite{A} et not\'es $D_{Q,\a}(y_1,\ldots,y_n)$ (il suffit, avec les notations de \cite{A}, 
de prendre $D=k(w_1,\ldots,w_r)$,  $Q=\L$ , et $\a_i$ est la restiction de $\sigma_i$ \`a $D$).

\medskip

Lorsque $n=2$, la matrice $\L$ est param\'etr\'ee par un seul scalaire $\l$.
On notera alors $T_{2,r}^{\l}(k)$ pour $T_{2,r}^{\L}(k)$.

Pour $n=r=1$,
on notera $U(k)$ l'alg\`ebre $T_{1,1}^{(1)}(k)$, engendr\'ee  par deux g\'en\'erateurs $w,y$ v\'erifiant $[w,y]=y$.
C'est l'alg\`ebre enveloppante de l'alg\`ebre de Lie r\'esoluble non ab\'elienne de dimension 2 sur $k$.
Il est clair que son corps de fractions $\Frac(U(k))$ est le corps de Weyl ${\mathcal D}_1(k)=\Frac(A_1(k))$.

\subsubsection{Ar\^etes eul\'eriennes.}

Afin de pouvoir repr\'esenter par un graphe les relations de commutation des 
alg\`ebres $T_{n,r}^{\L}(k)$ (et donc aussi des corps $\Frac(S_{n,r}^{\L}(k))$ avec les nouveaux g\'en\'erateurs $w_1,\ldots,w_r,y_1,\ldots,y_n$) on introduit une nouvelle convention  :

Une ar\^ete \'epaisse barr\'ee
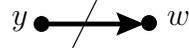
\begin{center}
\input{ar.pstex_t}
\end{center}
entre deux sommets $y$ et $w$ signifie que, dans l'alg\`ebre que l'on consid\`ere, on a entre les g\'en\'erateurs $y$ et $w$ la relation de commutation :
$[w,y]=y$, c'est-\`a-dire encore $yw=(w-1)y$.
Nous appellerons 
 ar\^etes  eul\'eriennes les ar\^etes de ce type (cette terminologie est motiv\'ee par le fait que dans un anneau de polyn\^omes c'est ce type de relations que v\'erifie une d\'erivation eul\'erienne par rapport \`a une variable avec l'op\'erateur de multiplication par cette variable). 





\subsubsection{Remarques.}

$\bullet$
De m\^eme que nous avons montr\'e ci-dessus au paragraphe \ref{redgrsec} que les alg\`ebres d\'efinies par g\'en\'erateurs et relations \`a partir d'un graphe $qW$, sont exactement les alg\`ebres polynomiales mixtes crois\'ees, se pose le probl\`eme de savoir quelles sont  les alg\`ebres correspondant  
 aux graphes compos\'es uniquement  d'ar\^etes quantiques et d'ar\^etes eul\'eriennes. Ce probl\`eme est r\'esolu au chapitre 7 de \cite{these} avec la d\'efinition d'alg\`ebre polynomiale ``mixte eul\'erienne''.

\medskip

$\bullet$
L'un des
int\'er\^ets de cette pr\'esentation des corps $\Frac(S^{\L}_{n,r}(k))$
par des graphes eul\'eriens est qu'elle permet de plonger les corps de fonctions rationnelles mixtes crois\'es 
 dans des corps de s\'eries de Malcev-Neuman, comme nous allons le voir au paragraphe suivant.

\subsection{Plongements dans des corps de s\'eries crois\'ees de Malcev-Neuman.} \label{plgsnrl}


{\bf Notations.}
En reprenant les notations du paragraphe \ref{arres}, on fixe un corps de fonctions rationnelles mixte crois\'e :
$$\Frac(S_{n,r}^{\L}(k))=k(w_1,\ldots,w_r)(y_1;\sigma_1)\ldots(y_n;\sigma_n).$$
On notera $k(w_1,\ldots,\what{w_i},\ldots,w_r)$ le sous-corps de $k(w_1,\ldots,w_r)$ engendr\'e sur $k$ par les variables $w_j$ pour $1\leq j\leq r$ et $j\neq i$.

On va construire un
 plongement de $\Frac(S_{n,r}^{\L}(k))$ 
 dans un corps 
 de s\'eries de Malcev-Neuman
adapt\'e, suivant en cela une m\'ethode d\'esormais classique (on pourra en particulier se r\'ef\'erer \`a l'article
 \cite{Ca1} de G. Cauchon, tr\`es complet et explicite sur la question;
 voir aussi \cite{AD3} et \cite{these}).


\subsubsection{Premi\`ere \'etape.} \label{blblbl}


Notons $s_i=w_i^{-1}$ pour tout $1\leq i\leq r$, de sorte que $D=k(w_1,\ldots,w_r)$ se plonge dans le corps de s\'eries de Malcev-Neuman commutatif $K=k_M[[s_1,\ldots,s_r]]$.
Les \'el\'ements de $K$ sont des s\'eries $T=\ssum_{a\in {\mathbb Z}^r} T(a)s^a$, avec $T(a)\in k$, dont le support est une partie bien ordonn\'ee de ${\mathbb Z}^r$. Pour un tel $T\in K^*$, on pose $\mu(T)={\rm min(supp}(T)) \in {\mathbb Z}^r$, et $\Psi(T)=T(\mu(T)) \in k^*$.
On v\'erifie facilement que :
\begin{center}
 $\mu(T_1T_2)=\mu(T_1)+\mu(T_2)$,\
et  \ $\Psi(T_1T_2)=\Psi(T_1)\Psi(T_2)$.
\end{center}
En particulier $\mu(T^{-1})=-\mu(T)$ et $\Psi(T^{-1})=\Psi(T)^{-1}$.
Pour tout $i$, $1\leq i\leq r$, on note $\mu_i(T)$ la $i^{\smmath{\`eme}}$ coordonn\'ee de $\mu(T)$.
Le $r$-uplet $\mu(T)$ est la valuation de $T$.
Le scalaire non-nul $\Psi(T)$ est le coefficient du terme de  valuation minimale de $T$.
Les restrictions \`a $D=k(w_1,\ldots,w_r)$ des  automorphismes $\sigma_1,\ldots,\sigma_r$  d\'efinis ci-dessus se prolongent en des automorphismes de $K$ encore not\'es $\sigma_1,\ldots,\sigma_r$  en posant $\sigma_i(s_i)=\ssum_{k\geq 1} s_i^k$ et $\sigma_i(s_j)=s_j$ si $j\neq i$. 

\begin{lemmeb}  
\begin{enumerate}
\item $\forall \ a=(a_1,\ldots,a_r)\in{\mathbb Z}^r,\ \forall \ T\in K^*,$\\
$ \Psi(\sigma_1^{a_1}\circ\ldots\sigma_r^{a_r}(T))=\Psi(T)$ ;
\item soit $T\in D^*$, et supposons qu'il existe $\a\in k^*$ tel que $\sigma_i(T)=\a T$. Alors $\a=1$ et $T$ appartient au corps  $k(w_1,\ldots,\what{w_i},\ldots,w_r)$.
\end{enumerate}
\end{lemmeb}
\begin{demo}
1. C'est clair par d\'efinition des $\sigma_i$.

2. Par le 1 on a $\Psi(\sigma_i(T))=\Psi(T)$ donc $\a=1$, et
 $\sigma_i(T)=T$. Dans l'extension $k(w_1,\ldots,\what{w_i},\ldots,w_r)((s_i))$ de $D$, on peut \'ecrire $T$ sous la forme $T=\ssum_{a\geq a_0}P(a)s_i^a$, avec $P(a)\in k(s_1,\ldots,\what{s_i},\ldots,s_r)$.
Alors : 
$$
\begin{array}{l}
T=P(a_0)s_i^{a_0}+P(a_0+1)s_i^{a_0+1}+\ldots,\\
\noalign{\medskip}
 \textrm{et }\ 
\sigma_i(T)=P(a_0)s_i^{a_0}+(P(a_0+1)+a_0P(a_0))s_i^{a_0+1}+\ldots.
\end{array}$$
 Donc $a_0=0$ et $T$ est de valuation nulle en $s_i$ dans $k(w_1,\ldots,\what{w_i},\ldots,w_r)((s_i))$. On recommence alors avec $T'=T-P(0)$, et n\'ecessairement $T'=0$, c'est-\`a-dire que $T=P(0)\in k(w_1,\ldots,\what{w_i},\ldots,w_r)$.
\end{demo}


\subsubsection{Deuxi\`eme \'etape.} \label{322}


On d\'efinit  un morphisme de groupes 

\begin{center}
$\sigma  :{\mathbb Z}^n \to \Aut(K)$ par $\sigma^{(a_1,\ldots,a_n)}=\sigma_1^{a_1}\circ\ldots\sigma_r^{a_r}$.
\end{center}

Notons que si $r<n$ les coefficients $a_{r+1},\ldots,a_n$ de $a$ ne jouent aucun r\^ole dans la d\'efinition ci-dessus.
L'application $\a  :{\mathbb Z}^n\times{\mathbb Z}^n \to k^*$ d\'efinie par 
$$\a(a,b)=\pprod_{1\leq i<j\leq n} \l_{i,j}^{a_ib_j}$$
 est un 2-cocycle. 

Le corps de fonctions rationnelles $\Frac(S^{\L}_{n,r}(k))$ se plonge alors dans le corps de s\'eries crois\'ees ${\mathbb F}=K_M[[{\mathbb Z}^n ;\sigma,\a]]$  d\'efini comme suit.
Les \'el\'ements de ${\mathbb F}$ sont des s\'eries $X=\ssum_{a\in {\mathbb Z}^n} X(a)y^a$, avec $X(a)\in K$, dont le support est une partie bien ordonn\'ee de ${\mathbb Z}^n$. Pour un tel $X\in {\mathbb F}^*$, on pose $\nu(X)={\rm min(supp}(X)) \in {\mathbb Z}^n$ la valuation de $X$, et $\Phi(T)=x(\nu(x)) \in K^*$ le coefficient de son terme de valuation minimale.
A nouveau on v\'erifie que :
\begin{center}
$\nu(T_1T_2)=\nu(T_1)+\nu(T_2)$,\
et \ $\Phi(T_1T_2)=\Phi(T_1)\sigma^{\nu(T_1)}(\Phi(T_2))\a(\nu(T_1),\nu(T_2))$.
\end{center}
En particulier $\nu(T^{-1})=-\nu(T)$\ et\ $\Phi(T^{-1})=\sigma^{-\nu(T)}(\Phi(T)^{-1})\a(\nu(T),\nu(T))$.

\medskip

Le plongement ci-dessus va nous permettre de d\'eterminer pour les corps de fonctions rationnelles mixtes crois\'es les premiers invariants rationnels  :
le centre, le groupe $G$ d\'efini dans \cite{AD1}, puis dans la section suivante un nouvel invariant  : le sous-tore quantique maximal simple.

\subsection{Centre des corps de fonctions rationnelles mixtes crois\'es.}

Pour une $k$-alg\`ebre $A$ on d\'esigne par ${\mathcal Z}(A)$ le centre de $A$.

\begin{spropo} \label{centre}
Soient $n\geq 1$ et $0\leq r\leq n$ deux entiers, et $\L\in M_n(k^*)$ multiplicativement antisym\'etrique.
 Le centre du corps $\Frac(S^{\L}_{n,r}(k))$ est l'intersection du centre du sous-corps engendr\'e par $y_1,\ldots,y_n$ avec le sous-corps engendr\'e par $y_{r+1},\ldots,y_n$  :
$${\mathcal Z}( \Frac(S^{\L}_{n,r}(k)))= {\mathcal Z}(k_{\L}(y_1,\ldots,y_n)) \cap k_{\wt\L}(y_{r+1},\ldots,y_n),$$
o\`u $\wt\L$ est la matrice extraite de $\L$ en \'eliminant les $r$ premi\`eres lignes et colonnes.
\end{spropo}
\begin{demo}
Soit $X=\ssum_{a\in{\mathbb Z}^n}X(a)y_1^{a_1}\ldots y_n^{a_n}$ un \'el\'ement du centre de $\Frac(S^{\L}_{n,r}(k))$, d\'evelopp\'e dans le sous-corps $D_M[[{\mathbb Z}^n;\sigma,\a]]$ de ${\mathbb F}$.
Pour tout $i$, $1\leq i\leq r$, on a  : $[w_i,X]=0$.
Or  :
$$[w_i,X]=\sum_aX(a)[w_i,y_1^{a_1}\ldots y_n^{a_n}]= \sum_aX(a)a_iy_1^{a_1}\ldots y_n^{a_n}.$$
Donc pour tout $i\leq r$, $X$ ne contient dans son d\'eveloppement aucun mon\^ome avec un exposant non-nul en $y_i$.
Ainsi
$X=\ssum_aX(a)y_{r+1}^{a_{r+1}}\ldots y_n^{a_n}$.

De m\^eme, pour tout $i$, $1\leq i\leq r$, on a  : $[y_i,X]=0$.
Or  :
$$[y_i,X]=\sum_a\left(\sigma_i(X(a))-\prod_{j\geq r+1}\l_{j,i}^{a_j}X(a)\right)y_i y_{r+1}^{a_{r+1}}\ldots y_n^{a_n}.$$
Donc   $\sigma_i(X(a))=\pprod_{j\geq r+1}\l_{j,i}^{a_j}X(a)$ pour tout $a\in {\mathbb Z}^n$, d'o\`u l'on d\'eduit  par le lemme \ref{blblbl} d'une part  que : 
$$X(a)\neq0 \oif \pprod_{j\geq r+1}\l_{j,i}^{a_j}=1,$$
 c'est-\`a-dire que $y_{r+1}^{a_{r+1}}\ldots y_n^{a_n}$ commute \`a $y_i$ d\`es que $X(a)\neq 0$, et d'autre part que :

\begin{center}
$X(a)\in k(w_1,\ldots,\what{w_i},\ldots,w_r)$ pour tout $a\in {\mathbb Z}^n$.
\end{center}

Puisque ceci est vrai pour tout $i\in \{1,\ldots,r\}$, on a $X(a)\in k$ pour tout $a\in{\mathbb Z}^n$.
Donc ${\mathcal Z}(\Frac(S^{\L}_{n,r}(k)))$ est inclus dans le corps de s\'eries $k_M[[y_{r+1},\ldots,y_n ;\a]]$. Par le th\'eor\`eme 1.10 de \cite{Ca1}, on en d\'eduit que les \'el\'ements centraux de $\Frac(S^{\L}_{n,r}(k))$ sont dans $k(y_{r+1},\ldots,y_n ;\a)$.
Par ailleurs un  \'el\'ement de ce centre 
doit commuter \`a tous les $y_i$, et donc appartenir \`a 
$${\mathcal Z}(k_{\L}(y_1,\ldots,y_n)) \cap k_{\wt\L}(y_{r+1},\ldots,y_n).$$

La r\'eciproque est \'evidente.
\end{demo}

\subsection{Invariants $E$ et $G$.}

Les  deux invariants suivants sont d\'efinis  dans \cite{AD1} afin de s\'eparer notamment les corps de Weyl des corps de fonctions rationnelles quantiques.
\begin{sdefin} \label{defeg}
Pour toute $k$-alg\`ebre $A$, on note :
\begin{itemize}
\item $G(A)=(A^*)'\cap k^*$ la trace sur $k^*$ du groupe d\'eriv\'e du groupe multiplicatif $A^*$ des unit\'es de $A$;
\item $E(A)=[A,A] \cap k$ la trace sur $k$ de l'alg\`ebre de Lie  d\'eriv\'ee   de $A$.
\end{itemize}
\end{sdefin}

{\bf Remarque.}
L'espace $E(A)$ ne peut que valoir $k$  ou $\{0\}$, suivant la pr\'esence ou pas dans  $A$ de deux \'el\'ements $x,y$ tels que $[x,y]=1$.

\begin{spropo} \label{eg}
Soient $n\geq 1$ et $0\leq r\leq n$ deux entiers, et $\L\in M_n(k^*)$ multiplicativement antisym\'etrique.
Alors le corps $\Dr$ v\'erifie les deux points suivants  :
\begin{enumerate}
\item $G(\Dr)=\ \langle \l_{i,j}\rangle $, le sous-groupe de $ k^*$ engendr\'e par les $\l_{i,j}$ ;
\item  
$E(\Dr)=
k {\textrm{ d\`es que }} r\geq 1$, et 
$E(\Frac(S_{n,0}^{\L}(k)))=0$.
\end{enumerate}
\end{spropo}
\begin{demo}
1. Il est clair que pour tout couple $(i,j)$ on a $y_iy_jy_i^{-1}y_j^{-1}=\l_{i,j}$, et donc $\langle \l_{i,j}\rangle $ est inclus dans $G(\Dr)$.

R\'eciproquement,  soit  $XYX^{-1}Y^{-1}$  un commutateur dans $\Dr$.
  On calcule 
$$\Phi(XYX^{-1}Y^{-1})= \Phi(XYX^{-1})\sigma^{\nu(XYX^{-1})}(\Phi(Y^{-1}))\a(\nu(XYX^{-1}),\nu(Y^{-1})).$$
 En r\'eit\'erant, et sachant que 
$$\Phi(X^{-1})=\a(\nu(X),\nu(X))\sigma^{-\nu(X)}(\Phi(X)^{-1}),$$
 on obtient :
\begin{equation} \label{hoho}
\Phi(XYX^{-1}Y^{-1})=\epsilon \Phi(X) \sigma^{\nu(X)} (\Phi(Y)) \sigma^{\nu(Y)} (\Phi(X)^{-1})\phi(Y)^{-1},
\end{equation}
avec
$$\begin{array}{l}
\epsilon=\a(\nu(Y),\nu(Y^{-1})) \a(\nu(XY),\nu(X^{-1}))\times \\
\blanc \a(\nu(X),\nu(Y)) \a(\nu(X),\nu(X)) \a(\nu(Y),\nu(Y)).
\end{array}$$
 On applique alors $\Psi$ \`a l'\'egalit\'e (\ref{hoho}). 
Rappelons  que pour tout $T\in K$ et tout $a\in{\mathbb Z}^n$ on a $\Psi(T^{-1})=\Psi(T)^{-1}$ et $\Psi(\sigma^a(T))=\Psi(T)$.
En remarquant que $\a$ est \`a valeurs dans $\langle \l_{i,j}\rangle $, on obtient que $\Psi\circ\Phi( XYX^{-1}Y^{-1})\in\ \langle \l_{i,j}\rangle $.
On conclut, en notant que $\Psi(\Phi(X_1X_2))=\a(\nu(X_1),\nu(X_2))\Psi(\Phi(X_1))\Psi(\Phi(X_2))$ pour tous $X_1,X_2\in {\mathbb F}^*$, que $G(\Dr)\subset \langle \l_{i,j}\rangle $.

2. Il est clair que si $r\geq 1$, alors $E(\Dr)=k$. Si $r=0$, alors $\Frac(S^{\L}_{n,0}(k))$ est un corps de fonctions rationnelles quantique, et le r\'esultat est \'enonc\'e au th\'eor\`eme 3.10 de \cite{AD1}.
\end{demo}

\medskip

{\bf Remarque.}
Avec la terminologie introduite en \ref{13}, le cas purement classique est celui o\`u l'invariant $G$ est trivial, et le cas purement quantique celui o\`u $E$ est trivial.

\section{Un invariant ``quantique'' : le sous-tore quantique maximal sim\-ple.}

\subsection{Propri\'et\'es pr\'eliminaires.}

La d\'efinition d'un tore quantique $k_{\L}[y_1^{\pm1},\ldots,y_n^{\pm1}]$ a \'et\'e rappel\'ee au point 3 de la d\'efinition \ref{dq}. Il s'agit d'une alg\`ebre de polyn\^omes de Laurent non-commutatifs, de dimension de Gelfand-Kirillov \'egale \`a $n$ (voir le point 5.1 de \cite{MCP}).
J.C. McConnell et  J.J. Pettit ont \'etabli dans \cite{MCP} les caract\'erisations sui\-vantes de la simplicit\'e  d'un tore quantique.

\begin{spropo}
    \label{simple}
Soient $n\geq 1$, et  $\L \in M_n(k^*)$ multiplicativement antisym\'etrique.
Pour le tore quantique ${\mathcal O}_{\L}({k^*}^n)$, les conditions sui\-vantes sont \'equi\-valentes :
\begin{enumerate}
\item ${\mathcal O}_{\L}({k^*}^n)$ est simple;

\item le centre de ${\mathcal O}_{\L}({k^*}^n)$ est r\'eduit \`a $k$;

\item il n'existe pas de $\a \in {\mathbb Z}^n \setminus \{ 0 \}$, tel que 

\centerline{$\forall j, \   1 \leq j \leq n, \ \l_{1,j}^{\a_1}\ldots \l_{n,j}^{\a_n} = 1$.}
\end{enumerate}
\end{spropo}
\begin{demo}
C'est la proposition 1.3 de \cite{MCP}.
\end{demo}

\medskip

Les propri\'et\'es suivantes concernant les morphismes de $k$-alg\`ebres entre tores quantiques apparaissent d\'ej\`a en partie dans \cite{MCP} et \cite{A2}, et figurent sous la forme suivante    dans  \cite{Retq}. 
Afin de fixer les notations nous les rappelons ici.

\medskip

{\bf Notation.}
Soient
  $\L \in M_n(k^*)$ et $ \L' \in M_{n'}(k^*)$ multiplicativement antisym\'etriques, et les tores quantiques associ\'es
 ${\mathcal O}_{\L}({k^*}^n)=k_{\L}[y_1^{\pm 1},\ldots,y_n^{\pm 1}]$
et ${\mathcal O}_{\L'}({k^*}^{n'})=k_{\L'}[{y'_1}^{\pm 1},\ldots,{y'}^{\pm 1}_{n'}].$
 Soit un
 morphisme $\Phi  :{\mathcal O}_{\L}({k^*}^n) \to {\mathcal O}_{\L'}({k^*}^{n'})$.
L'\'el\'ement $\Phi (y_i)$  est inversible dans ${\mathcal O}_{\L'}({k^*}^{n'})$ pour tout $i \in \{1,\ldots,n\}$, c'est donc un mon\^ome $\a_i ({y'_1})^{h_{1,i}} \ldots ({y'_{n'}})^{h_{n',i}}$, o\`u $\a_i \in k^*$, et $h_{k,i} \in {\mathbb Z}$. On note alors $H_{y,y'}(\Phi)$ la matrice  $(h_{i,j}) \in M_{n',n}({\mathbb Z})$.

\begin{spropo}
\label{MP}
Soient $n,n',n''\geq 1$ trois entiers.
\begin{enumerate}
\item
 Soient $\L \in M_n(k^*)$ et $\L' \in M_{n'}(k^*)$ deux matrices multiplicativement antisym\'etriques, et les tores quantiques associ\'es
$${\mathcal O}_{\L}({k^*}^n)=k_{\L}[y_1^{\pm 1},\ldots,y_n^{\pm 1}] \ \textrm{ et }\ {\mathcal O}_{\L'}({k^*}^{n'}) =k_{\L'}[{y'_1}^{\pm 1},\ldots,{y'_{n'}}^{\pm 1}].$$
Soit $\Phi  : {\mathcal O}_{\L}({k^*}^n) \to {\mathcal O}_{\L'}({k^*}^{n'})$ un  morphisme. Notons $H_{y,y'}({\Phi}) = (h_{i,j})$. Alors on a :  
\begin{equation} \label{eqno1}
\forall i,j, \ 1\leq i,j \leq n, \   \ \ \ \l_{i,j} = \prod_{1 \leq k,t \leq n'} {{\l}'_{k,t}}^{h_{k,i}h_{t,j}}.  
\end{equation}
R\'eciproquement, \'etant donn\'ee une matrice $H= (h_{i,j}) \in M_{n',n}({\mathbb Z})$ satisfaisant les \'equations (\ref{eqno1}), il existe un unique morphisme  $\Phi$ de ${\mathcal O}_{\L}({k^*}^n)$ dans  ${\mathcal O}_{\L'}({k^*}^{n'})$ tel que $\Phi (y_i)={y'_1}^{h_{1,i}} \ldots {y'_n}^{h_{n,i}}$. Il v\'erifie  $H = H_{y,y'}({\Phi})$.

\item Soient $\L \in M_n(k^*),\ \L' \in M_{n'}(k^*)$, et $ \L'' \in M_{n''}(k^*)$ trois matrices multiplicativement antisym\'etriques, et les trois tores quantiques associ\'es 
$$
\begin{array}{l}
\tl({k^*}^n)=k_{\L}[y_1^{\pm 1},\ldots,y_n^{\pm 1}],
\quad  {\mathcal O}_{\L'}({k^*}^{n'})=k_{\L'}[{y'_1}^{\pm 1},\ldots,{y'}^{\pm 1}_{n'}],\\
\noalign{\medskip}
\textrm{et }\ {\mathcal O}_{\L''}({k^*}^{n''})=k_{\L''}[{y''_1}^{\pm 1},\ldots,{y''}^{\pm 1}_{n''}].
\end{array}$$
 Consid\'erons en outre  deux  morphismes de $k$-alg\`ebres 
$$\Phi  : \tl({k^*}^n) \to{\mathcal O}_{\L'}({k^*}^{n'}) \ \textrm{ et } \ \Phi'  : {\mathcal O}_{\L'}({k^*}^{n'}) \to {\mathcal O}_{\L''}({k^*}^{n''}).$$
 Alors :
$$H_{y,y''}({\Phi' \circ \Phi}) = H_{y',y''}({\Phi'}) H_{y,y'}({\Phi}),$$
 au sens du produit usuel des matrices.

\item Soient $\L \in M_n(k^*)$ et $ \L' \in M_{n'}(k^*)$ deux matrices multiplicativement antisym\'etriques, et les tores quantiques associ\'es 
$$\tl({k^*}^n)=k_{\L}[y_1^{\pm 1},\ldots,y_n^{\pm 1}]\ \textrm{ et } \  {\mathcal O}_{\L'}({k^*}^{n'})=k_{\L'}[{y'_1}^{\pm 1},\ldots,{y'_{n'}}^{\pm 1}].$$
 Soit $\Phi  : \tl({k^*}^n) \to{\mathcal O}_{\L'}({k^*}^{n'})$ un  morphisme. Alors $\Phi$ est un isomorphisme si et seulement si $n=n'$ et $H_{y,y'}({\Phi}) \in GL_n({\mathbb Z})$.

\item  En particulier, deux tores quantiques $\tl({k^*}^n)$ et ${\mathcal O}_{\L'}({k^*}^{n})$ sont $k$-i\-so\-mor\-phes si et seulement s'il existe une matrice $H \in GL_n({\mathbb Z})$ v\'erifiant les relations (\ref{eqno1}).
\end{enumerate}
\end{spropo}
\begin{demo}
Il s'agit des points (ii) \`a (v) du lemme 1.4 de \cite{Retq}.
\end{demo}

\medskip

Par ailleurs on d\'emontre dans \cite{Retq} le th\'eor\`eme suivant.

\begin{stheoreme} \label{thpr}
 Soient $n\geq 1$ un entier, et $\L \in M_n(k^*)$ multiplicativement antisym\'etrique. Supposons le tore quantique  ${\mathcal O}_{\L}({k^*}^n)$  simple. Alors tout endomorphisme de ${\mathcal O}_{\L}({k^*}^n)$ est un automorphisme de ${\mathcal O}_{\L}({k^*}^n)$.
\end{stheoreme}
\begin{demo}
C'est le th\'eor\`eme 3.6 de \cite{Retq}.
\end{demo}

\subsection{Sous-tores quantiques maximaux.}

\begin{sdefin} \label{defpp}
Soient $n\geq 1$ un entier, et $\L\in M_n(k^*)$ une matrice multiplicativement antisym\'etrique.
Soit $A$ une $k$-alg\`ebre.
On dit que ${\mathcal O}_{\L}({k^*}^n)$ est un sous-tore quantique maximal de $A$ si les deux conditions suivantes sont r\'ealis\'ees  :
\begin{enumerate}
\item il existe une sous-alg\`ebre de $A$  $k$-isomorphe \`a ${\mathcal O}_{\L}({k^*}^n)$;
\item
 pour tout tore quantique  ${\mathcal O}_M({k^*}^m)$, avec $m\geq 1$ un entier et $M\in M_m(k^*)$ multiplicativement antisym\'etrique, 
 s'il existe un morphisme d'alg\`ebres 
$$\phi :{\mathcal O}_M({k^*}^m) \rightarrow A,$$
 alors il existe un morphisme $\wt\phi :{\mathcal O}_M({k^*}^m) \rightarrow {\mathcal O}_{\L}({k^*}^n)$.
\end{enumerate}
\end{sdefin}

Avec ces notations,
le point 2 signifie  que, si ${\mathcal O}_{\L}({k^*}^n)$ est un sous-tore quantique maximal de $A$, alors tout tore quantique admettant une image homomorphe dans $A$ admet une image homomorphe dans ${\mathcal O}_{\L}({k^*}^n)$.
En particulier, si ${\mathcal O}_{\L}({k^*}^n)$ est un sous-tore quantique maximal de $A$, alors tout tore quantique simple qui se plonge dans $A$ se plonge dans ${\mathcal O}_{\L}({k^*}^n)$.

\medskip

{\bf Remarques.}
$\bullet$ Un tore quantique est toujours sous-tore quantique maximal de lui-m\^eme.

$\bullet$
Un sous-tore quantique maximal n'est pas unique :
par exemple, il est facile \`a partir de la proposition \ref{MP} de v\'erifier que le tore quantique multiparam\'etr\'e ${\mathcal O}_{\L}({k^*}^3)=k_{\L}[y_1^{\pm 1};y_2^{\pm 1};y_3^{\pm 1}]$, avec $\L=\left({\tiny \begin{array}{ccc}
1&q&1\\
q^{-1}&1&1\\
1&1&1 \end{array}}\right)$,
admet au moins deux sous-tores quantiques maximaux :
lui-m\^eme et le tore quantique $k_q[y_1^{\pm 1},y_2^{\pm 1}]$.
Par contre, le lemme suivant, bas\'e sur le th\'eor\`eme  \ref{thpr},  permet de d\'emontrer qu'un sous-tore quantique maximal {\sl simple}, quand il existe, est unique \`a $k$-isomorphisme pr\`es.

\begin{slemme} \label{kifoqr}
Soit $n\geq 1$ un entier.
Soient $\L,\L'\in M_n(k^*)$ multiplicativement antisym\'etriques.
Soient $A$ et $B$ deux $k$-alg\`ebres.
Supposons que ${\mathcal O}_{\L}({k^*}^n)$ est un sous-tore quantique maximal de
$A$, et ${\mathcal O}_{\L'}({k^*}^n)$ un sous-tore quantique maximal de $B$.
Si $A$ et $B$ sont $k$-isomorphes, et si ${\mathcal O}_{\L}({k^*}^n)$ est simple, alors ${\mathcal O}_{\L}({k^*}^n)$ et ${\mathcal O}_{\L'}({k^*}^n)$ sont $k$-isomorphes.
\end{slemme}
\begin{demo}
L'isomorphisme de $A$ sur $B$ induit un morphisme injectif  $\phi$ de ${\mathcal O}_{\L}({k^*}^n)$ dans $B$, donc un morphisme  $\wt\phi$ de ${\mathcal O}_{\L}({k^*}^n)$  dans ${\mathcal O}_{\L'}({k^*}^n)$. De m\^eme on a un morphisme $\wt\psi$ de ${\mathcal O}_{\L'}({k^*}^n)$ dans ${\mathcal O}_{\L}({k^*}^n)$.
Donc $\wt\psi \circ\wt\phi$ est un endomorphisme de ${\mathcal O}_{\L}({k^*}^n)$.
Puisque celui-ci est suppos\'e simple, $\wt\psi \circ\wt\phi$ est un automorphisme de ${\mathcal O}_{\L}({k^*}^n)$ d'apr\`es le th\'eor\`eme \ref{thpr}.
D'apr\`es la proposition \ref{MP} les matrices $H$ et $H'$ de $M_n({\mathbb Z})$ associ\'ees \`a $\wt\psi$ et $\wt\phi$ ont donc leur produit dans $GL_n({\mathbb Z})$, ce qui implique que $H\in GL_n({\mathbb  Z})$, et donc, toujours d'apr\`es la proposition \ref{MP}, $\wt\psi$ est un isomorphisme.
\end{demo}

\begin{spropo} \label{prehappyqr}
Soient $n,n'\geq 1$ des entiers, et  $\L\in M_n(k^*)$ et $\L'\in M_{n'}(k^*)$ deux matrices multiplicativement antisym\'etriques.
Soient $A$ et $B$ deux $k$-alg\`ebres.
Si ${\mathcal O}_{\L}({k^*}^n)$ est un sous-tore quantique maximal \emph{simple} de $A$, si ${\mathcal O}_{\L'}({k^*}^{n'})$ est un sous-tore quantique maximal \emph{simple} de $B$, et si $A$ et $B$ sont $k$-isomorphes, alors
$n=n'$ et les tores quantiques  ${\mathcal O}_{\L}({k^*}^n)$ et ${\mathcal O}_{\L'}({k^*}^n)$ sont $k$-isomorphes.
\end{spropo}
\begin{demo}
Comme dans la preuve du  lemme \ref{kifoqr},  l'isomorphisme entre $A$ et $B$ induit deux morphismes $\wt\phi$ de ${\mathcal O}_{\L}({k^*}^n)$  dans ${\mathcal O}_{\L'}({k^*}^{n'})$ et $\wt\psi$ de ${\mathcal O}_{\L'}({k^*}^{n'})$ dans ${\mathcal O}_{\L}({k^*}^n)$.
Mais puisque ces tores quantiques sont simples, les deux morphismes consid\'er\'es sont injectifs, d'o\`u il d\'ecoule que  ces deux tores quantiques ont  la m\^eme GK-dimension, c'est-\`a-dire que $n=n'$. 
On conclut alors \`a l'aide du lemme \ref{kifoqr}.
\end{demo}

\begin{scorol}
Si une $k$-alg\`ebre admet un sous-tore quantique maximal \emph{sim\-ple}, alors ce dernier est unique \`a isomorphisme pr\`es.
\end{scorol}
\begin{demo}
On applique la proposition pr\'ec\'edente avec $A=B$.
\end{demo}

\medskip

A la suite de la remarque faite apr\`es la d\'efinition \ref{defpp}, notons qu'un tore quantique simple est toujours un sous-tore quantique maximal simple de lui-m\^eme.
Rappelons \'egalement qu'un isomorphisme entre deux tores quantiques  s'interpr\`ete gr\^ace \`a la proposition \ref{MP} comme une condition sur les coefficients des matrices de param\'etrisation.

\medskip

{\bf Remarque.}
Comme le montre la proposition suivante, une $k$-alg\`ebre n'admet pas toujours un sous-tore quantique maximal simple.

\begin{spropo}[Contre-exemple]
Soit $q\in k^*$ non racine de l'unit\'e.
Soit  $T={\mathcal O}_{\L}({k^*}^4)$, param\'etr\'e par 
$$\L=\left(\begin{array}{cccc}
1&q&1&1\\
q^{-1}&1&1&1\\
1&1&1&-1\\
1&1&-1&1 \end{array}\right).$$
Alors $T$ n'admet pas de sous-tore quantique maximal simple.
\end{spropo}
\begin{demo}
Supposons que $D$ admette un sous-tore quantique maximal simple $S$.
D'apr\`es la remarque suivant la d\'efinition \ref{defpp}, $S$ doit contenir les tores quantiques simples ${\mathcal O}_q({k^*}^2)$ et  ${\mathcal O}_{-q}({k^*}^2)$.
Donc le groupe $G(S)$ doit contenir $q$ et $-q$.
Mais $G(S)\subset G(T)=\ \langle q,-1\rangle $, donc $G(S)=G(T)$.
Ainsi $S$ n'est pas uniparam\'etr\'e, et donc $\gkd(S)\geq 3$.

Supposons que $\gkd(S)=3$.
Notons $S=k_M[y_1^{\pm 1},y_2^{\pm 1},y_3^{\pm 1}]$, param\'etr\'e par  $M=(\mu_{i,j})\in M_3(k^*)$ multiplicativement antisym\'etrique.
Consid\'erons la sous-alg\`ebre 
 $S'={k}_{M'}[y_1^{\pm 2},y_2^{\pm 2},y_3^{\pm 2}]$
 engendr\'ee par les carr\'es des g\'en\'erateurs de $S$, o\`u $M'=(\mu_{i,j}^4)$.
D'apr\`es la proposition \ref{simple} appliqu\'ee aux coefficients de $M$ et de $M'$, on peut v\'erifier  que la simplicit\'e de  $S$ implique celle de $S'$.
Or de  $G(S)=\ \langle -1,q\rangle $ il d\'ecoule que $G(S')=\langle q^4\rangle $. Le tore quantique $S'$ serait   uniparam\'etr\'e et de GK-dimension \'egale \`a 3, et  ne pourrait pas \^etre simple en vertu de la proposition 2.3(ii) de \cite{Retq}.
On aboutit donc \`a une contradiction en supposant $S$ de GK-dimension \'egale \`a 3.

\smallskip

Supposons que $\gkd(S)=4$.
Comme pr\'ec\'edemment notons 
$$S=k_M[y_1^{\pm 1},y_2^{\pm 1},y_3^{\pm 1},y_4^{\pm 1}],$$
 avec $M=(\mu_{i,j})\in M_4(k^*)$ multiplicativement antisym\'etrique.
La sous-alg\`ebre 
 $S'={k}_{M'}[y_1^{\pm 2},y_2^{\pm 2},y_3^{\pm 2},y_4^{\pm 2}]$ engendr\'ee par les carr\'es des g\'en\'erateurs de $S$, o\`u $M'=(\mu_{i,j}^4)$, est un tore quantique simple uniparam\'etr\'e par $q^4$.
D'apr\`es la preuve de la proposition 2.3 de \cite{Retq}, $S'$ est $k$-isomorphe \`a un tore quantique simple ${\mathcal O}_{\L''}({k^*}^4)$ o\`u $\L''$ est de la forme canonique
$$\L''=\left({\scriptsize \begin{array}{cccc}
1&q^{4d_1}&1&1\\
q^{-4d_1}&1&1&1\\
1&1&1&q^{4d_2}\\
1&1&q^{-4d_2}&1 \end{array}}\right),$$
avec des entiers $d_1,d_2\geq 1$ tels que $d_1$ divise $d_2$.
Remarquons que les g\'en\'erateurs de $S$ sont inversibles dans $T$, ce sont donc des mon\^omes en les g\'en\'erateurs de $T$.
Par cons\'equent $S'$ est une sous-alg\`ebre du tore quantique $T'={\mathcal O}_{\L'}({k^*}^4)$ engendr\'e par les carr\'es des g\'en\'erateurs de $T$, avec 
$$\L'=\left({\scriptsize \begin{array}{cccc}
1&q^4&1&1\\
q^{-4}&1&1&1\\
1&1&1&1\\
1&1&1&1 \end{array}}\right).$$
Or $\Frac({\mathcal O}_{\L''}({k^*}^4))$ ne peut pas se plonger dans $\Frac({\mathcal O}_{\L'}({k^*}^4))$ d'apr\`es le corollaire 2.14 de \cite{AD3}.
A nouveau on aboutit \`a une contradiction.

Puisque $\gkd(T)=4$ on a ainsi d\'emontr\'e la proposition.
\end{demo}

\subsection{Cas des corps de fonctions rationnelles mixtes crois\'es.}

\begin{spropo}    \label{happy}
Soient $n\geq 1$ et $0\leq r\leq n$ deux  entiers, et $\L\in M_n(k^*)$ multiplicativement antisym\'etrique. 
Alors le tore quantique ${\mathcal O}_{\L}({k^*}^n)$ est un sous-tore quantique maximal de $\Dr$.
\end{spropo}
\begin{demo}
Il est clair que les $y_i$ et leurs inverses engendrent dans $\Dr$ un tore quantique isomorphe \`a ${\mathcal O}_{\L}({k^*}^n)$.
Consid\'erons alors un tore quantique  ${\mathcal O}_M({k^*}^m)$, avec $M=(\mu_{i,j})\in M_m(k^*)$ multiplicativement antisym\'etrique, et supposons qu'il existe  un morphisme $\phi$ de ${\mathcal O}_M({k^*}^m)$ dans $\Dr$. Notons $X_1^{\pm 1},\ldots,X_m^{\pm1}$ les images par $\phi$ des g\'en\'erateurs de ${\mathcal O}_M({k^*}^m)$. 
On a alors dans le corps $\Dr$ les identit\'es $X_iX_j=\mu_{i,j}X_jX_i$ pour tous $1\leq i,j\leq m$, et en appliquant dans le corps de s\'eries ${\mathbb F}$ d\'efini en \ref{322}  l'application $\Psi\circ\Phi$, on obtient 
$$\a(\nu(X_i),\nu(X_j))=\mu_{i,j}\a(\nu(X_j),\nu(X_i)).$$
Ainsi :
$$\forall \ 1\leq i,j\leq m, \ \mu_{i,j}=\prod_{1\leq k,t\leq n}\l_{k,t}^{\nu_k(X_i)\nu_t(X_j)}.$$
En notant 
$$H=(\nu_k(X_i))_{1\leq k\leq n, \  1\leq i\leq m},$$
 on en d\'eduit par la proposition \ref{MP} qu'il existe un morphisme de $k$-alg\`ebres de ${\mathcal O}_M({k^*}^m)$ dans ${\mathcal O}_{\L}({k^*}^n)$. 
\end{demo}

\begin{scorol}
Reprenons les hypoth\`eses de la proposition pr\'ec\'edente.
Supposons de plus que la matrice $\L$ v\'erifie l'une des trois conditions \'equivalentes de la proposition \ref{simple}.
 Alors
  ${\mathcal O}_{\L}({k^*}^n)$ est  un sous-tore quantique maximal simple de   $\Dr$.
\end{scorol}
\begin{demo}
Ceci est une cons\'equence directe des   propositions \ref{happy} et \ref{simple}.
\end{demo}

\section{Un invariant ``classique'' : le w-degr\'e sup\'erieur.} \label{secwdeg}

\subsection{Notion de corps de Weyl mixtes.}

On introduit ici  une version multiparam\'etr\'ee des corps de Weyl mixtes d\'efinis dans \cite{AD3}.

\begin{sdefin} \label{defwm}
Soient $m,n,t \in {\mathbb N}$, et soit $\qb=(q_1,\ldots,q_n)\in (k\setminus\{0,1\})^n$. 
On appelle corps de Weyl mixte associ\'e \`a ces donn\'ees, et on note ${\mathcal D}_{m,n,t}^{\qb}(k)$, le corps de fractions de l'alg\`ebre ${\mathcal A}_{m,n,t}^{\qb}(k)$, engendr\'ee sur $k$ par  $2m+2n+t$ g\'en\'erateurs 
$$x_1,\ldots,x_m,y_1,\ldots,y_m,u_1,\ldots,u_n,v_1,\ldots,v_n,z_1,\ldots,z_t$$
 soumis aux relations suivantes pour tous $1\leq i\neq j\leq m$, $1\leq k\neq l\leq n$, $1\leq p\neq s\leq t$ :
\begin{eqnarray*} 
&  [x_i,x_j]=[y_i,y_j]=[x_i,y_j]=0,& \\ 
 & x_iy_i=y_ix_i+1,& \\
 & {[u_k,u_l]}=[v_k,v_l]=[u_k,v_l]=0,& \\
 & u_kv_k=q_kv_ku_k,& \\
 & {[z_p,z_s]}=[z_p,x_i]=[z_p,y_i]=[z_p,u_k]=[z_p,v_k]=0,& \\ 
 & {[x_i,u_k]}=[x_i,v_k]=[y_i,u_k]=[y_i,v_k]=0.&
\end{eqnarray*}
\end{sdefin}

{\bf Remarque.}
Ces relations sont donc 
 repr\'esent\'ees par le graphe suivant :

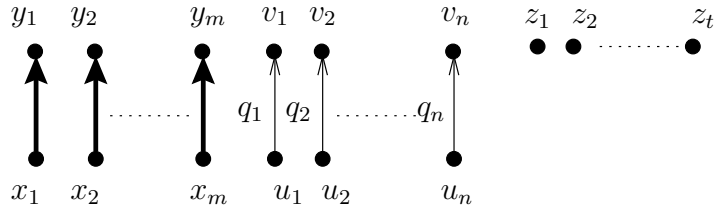
\begin{center}
\input{cwm.pstex_t}
\end{center}

{\bf Exemples.}
$\bullet$ Il est clair qu'on retrouve dans le cas $n=0$ les corps de Weyl classiques  :
${\mathcal D}_{m,0,t}(k)={\mathcal D}_{m,t}(k)=\Frac(A_{m,t}(k))$.

\smallskip

$\bullet$ Si $m=0$,
on a  :  ${\mathcal D}_{0,n,t}^{\qb}(k)=\Frac({\mathcal O}_{\L}(k^{2n+t}))$,  corps de fonctions rationnelles quantique param\'etr\'e par la matrice 
 $\L\in M_{2n+t}(k^*)$, compos\'ee sur sa diagonale  de $n$ blocs  ${\tiny \left(\begin{array}{cc} 1&q_i\\ q_i^{-1}&1 \end{array} \right)}$, et de 1 partout ailleurs.
Rappelons en particulier que, d'apr\`es le th\'eor\`eme 2.19 de \cite{Pa1}, tout corps de fractions d'un  tore quantique uniparam\'etr\'e est de ce type.



\begin{spropo} \label{heyhey}
Soient $m,n,t,\qb$ tels que dans la d\'efinition \ref{defwm}.
Alors :
$$\gkdpr({\mathcal A}_{m,n,t}^{\qb}(k))=\gktdpr({\mathcal D}_{m,n,t}^{\qb}(k))=2m+t+2n.$$
\end{spropo}
\begin{demo}
Comme \`a la proposition \ref{gkdsnrl}, on utilise r\'ecursivement le lemme 2.2 de \cite{HK}, puis le th\'eor\`eme 7.3 de \cite{Z}.
\end{demo}

\begin{spropo} \label{centrewm}
Soient $m,n,t,\qb$ tels que dans la d\'efinition \ref{defwm}.
Soit $r$ le nombre de $q_i$ racines de l'unit\'e dans $k^*$.
Alors ${\mathcal Z}({\mathcal D}_{m,n,t}^{\qb}(k))$ est une extension transcendante pure de $k$, de degr\'e $r+t$.

En particulier ${\mathcal Z}({\mathcal D}_{m,n,t}^{\qb}(k))=k(z_1,\ldots,z_t)$ d\`es que tous les $q_i$ sont non racines de l'unit\'e.
\end{spropo}
\begin{demo}
C'est un cas particulier de la proposition \ref{centre}.
\end{demo}

\subsection{Un invariant dimensionnel : le w-degr\'e inf\'erieur.}

{\bf Exemple pr\'eliminaire.}
Fixons $q\in k^*$ non racine de l'unit\'e, et consid\'erons les deux corps de Weyl mixtes $F={\mathcal D}_{1,0,2}^{(q,q)}(k)$ et $F'={\mathcal D}_{2,0,1}^{(q)}(k)$ :
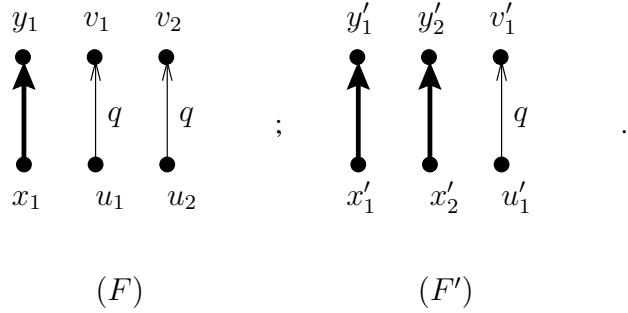
\begin{center}
\input{expre.pstex_t}
\end{center}

Il est clair, au vu de tous les r\'esultats qui pr\'ec\`edent,
que 
\begin{itemize}
\item $\gktd(F)=\gktd(F')=6$;
\item ${\mathcal Z}(F)={\mathcal Z}(F')=k$;
\item $E(F)=E(F')=k$.
\end{itemize}
De plus, d'apr\`es la proposition  
 \ref{eg}  on a $G(F)=G(F')=\langle q\rangle$ (r\'esultat d\'ej\`a montr\'e \`a la proposition 2.5 de \cite{AD3}).

L'invariant suivant, d\'efini dans le cas uniparam\'etr\'e dans l'article \cite{AD3}, a entre autres pour objectif de s\'eparer ce type de situations.

\begin{sdefin}[w-degr\'e inf\'erieur]
Soit $L$ un corps gauche  contenant $k$ dans son centre.
Notons $M$ la borne sup\'erieure de l'ensemble des entiers $m\geq 1$ pour lesquels il existe dans $L$ un sous-corps $k$-isomorphe au corps de Weyl classique ${\mathcal D}_{m}(k)$, avec la convention $M=0$ s'il n'existe pas de tels sous-corps.
Alors $2M$ est  appel\'e {\rm w}-degr\'e inf\'erieur de $L$, et  not\'e $\winfpr(L)$.
\end{sdefin}

{\bf Remarques.} $\bullet$
En d'autres termes, le w-degr\'e inf\'erieur de $L$ est  le GK-degr\'e de transcendance  du plus grand corps de Weyl classique  qui se plonge dans $L$.
On a en particulier $\winf(L)\leq \gktd(L)$.

\smallskip

$\bullet$
J. Alev et F. Dumas calculent dans \cite{AD3} le w-degr\'e inf\'erieur d'un corps de Weyl mixte dans le cas o\`u les $q_i$ sont tous puissances d'un m\^eme $q\in k^*$ non racine de l'unit\'e.
Leur m\'ethode de preuve (par r\'eduction modulo $p$) s'adapte sans difficult\'e au cas multiparam\'etr\'e.
On obtient ainsi le th\'eor\`eme suivant.

\begin{stheoreme} \label{thwdeg}
Soient $m,n,t,\qb$ tels que dans la d\'efinition \ref{defwm}.
Alors : 
$$\winfpr({\mathcal D}_{m,n,t}^{\qb}(k))=2m.$$
\end{stheoreme}
\begin{demo}
On adapte la preuve du th\'eor\`eme 2.11 de \cite{AD3}. Pour plus de d\'etails, voir la preuve du  th\'eor\`eme 3.1.2.4 de \cite{these}.
\end{demo}

\medskip

Ainsi le w-degr\'e inf\'erieur ``caract\'erise'' le nombre d'ar\^etes de Weyl d'un corps de Weyl mixte, et deux corps de Weyl mixtes $k$-isomorphes doivent avoir le m\^eme nombre d'ar\^etes de Weyl.
Ceci r\'esoud en particulier l'exemple pr\'esent\'e au d\'ebut de ce paragraphe, en montrant que les corps $F$ et $F'$ ne sont pas $k$-isomorphes.

A titre d'illustration, on d\'eduit du th\'eor\`eme \ref{thwdeg} le r\'esultat de plongement suivant (valable en particulier pour les corps de Weyl classiques).

\begin{scorol}
Avec les notations de la d\'efinition \ref{defwm},
si un corps de Weyl mixte ${\mathcal D}_{m,n,t}^{\qb}(k)$ se plonge dans un corps de Weyl mixte ${\mathcal D}_{m',n',t'}^{\surl q'}(k)$, alors $m\leq m'$.
\end{scorol}
\begin{demo}
Ceci d\'ecoule directement de la d\'efinition du w-degr\'e inf\'erieur et du th\'eor\`eme \ref{thwdeg}.
\end{demo}

\subsection{Notion de w-degr\'e sup\'erieur.}

Nous introduisons dans ce paragraphe un autre invariant dimensionnel, directement li\'e au w-degr\'e inf\'erieur, mais plus facile \`a calculer que ce dernier pour les corps de fonctions rationnelles mixtes crois\'es g\'en\'eraux.
On commence par un exemple.

\subsubsection{Un exemple de ``d\'etressage par plongement''.}\label{plg4}


Fixons  $q\in k^*$ non racine de l'unit\'e.
Consid\'erons l'alg\`ebre   $S_{2,2}^{\L}(k)$, avec $\L=\left({\tiny \begin{array}{cc}
1& q\\
q^{-1}&1
\end{array}}\right)$. 
Il s'agit  de l'alg\`ebre, not\'ee $S_{2,2}^q(k)$ avec les conventions introduites au paragraphe \ref{defqw}, engendr\'ee sur $k$ par $x_1,x_2,y_1,y_2$ avec les relations :

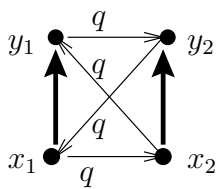
\begin{center}
\input{graf1.pstex_t}
\end{center}

Par des consid\'erations tr\`es techniques sur les centralisateurs de certaines paires d'\'el\'ements, il est d\'emontr\'e au th\'eor\`eme 3.5 de \cite{AD3} que $\Frac(S_{2,2}^{q}(k))$ n'est pas isomorphe \`a un corps de Weyl mixte (nous donnons un peu plus loin en \ref{eqwmsc4} une autre preuve de ce r\'esultat).
En revanche, il se plonge dans un corps de Weyl mixte de la fa\c con suivante.

\begin{propb} 
Il existe un plongement de $\Frac(S_{2,2}^{q}(k))$  dans le corps de Weyl mixte ${\mathcal D}_{2,1,0}^{q}(k)$.
\end{propb}
\begin{demo}
Il est facile de v\'erifier que le graphe suivant est admissible, et d\'efinit donc une extension it\'er\'ee de Ore $R$ engendr\'ee sur $k$ par 6 g\'en\'erateurs $x_1,x_2,y_1,y_2,y_3,y_4$ avec les relations de commutation :

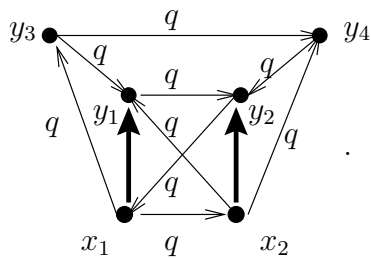
\begin{center}
\input{graf2.pstex_t}
\end{center}

Il est clair que l'alg\`ebre $S_{2,2}^{q}(k)$ se plonge dans $R$.
Mais $\Frac(R)$ admet comme   autre famille de g\'en\'erateurs :

\begin{center}
\input{graf3.pstex_t}
\end{center}
d'o\`u le r\'esultat annonc\'e.
\end{demo}


\subsubsection{D\'efinition  du w-degr\'e sup\'erieur.}\label{wwcwm}


 Le w-degr\'e inf\'erieur s'av\`ere \^etre un invariant difficile \`a calculer, et hormis dans le cas des corps de Weyl mixtes (th\'eor\`eme \ref{thwdeg}) et de certains exemples particuliers (voir la question (ii) \`a la fin  de \cite{AD3}), on ne dispose \`a notre connaissance d'aucune m\'ethode permettant son calcul.
Reprenons par exemple le corps $F=\Frac(S_{2,2}^{q}(k))$ ci-dessus.
Bien que l'on puisse raisonnablement conjecturer  que $\winf(F)=2$, la question d\'ej\`a pos\'ee en \cite{AD3} de d\'emontrer ce r\'esultat reste \`a notre connaissance toujours ouverte.
 Par contre, on a su en \ref{plg4} plonger $F$ dans un corps de Weyl mixte 
sans augmenter le nombre d'ar\^etes de Weyl dans le graphe traduisant les relations de commutation entre les g\'en\'erateurs du  corps.

\medskip

L'id\'ee est alors la suivante  :
puisqu'on ne sait pas mesurer ``par en bas'' le corps de Weyl classique maximal qu'on peut plonger dans le corps $F$, on va le faire ``par en haut'' en  plongeant $F$ dans un corps de Weyl mixte avec un nombre minimal d'ar\^etes de Weyl.
D'o\`u la notion suivante :

\begin{definb}
[w-degr\'e sup\'erieur]
Soit $L$ un corps gauche  contenant $k$ dans son centre.
Notons $M$ le plus petit  des entiers $m\geq 0$ pour lesquels il existe un corps de Weyl mixte dont le \emph{w}-degr\'e inf\'erieur est  \'egal \`a $2m$, et admettant un sous-corps $k$-isomorphe \`a $L$, avec la convention $M=+\infty$ s'il n'existe pas de tels corps de Weyl mixtes.
Alors $2M$ est  appel\'e {\rm w}-degr\'e sup\'erieur de $L$, et  not\'e $\wsuppr(L)$.
\end{definb}

Remarquons que, \`a cause du th\'eor\`eme \ref{thwdeg} et par d\'efinitions, on a toujours
$$
\winf(L)\leq\wsup(L).$$

\begin{propb} 
Soient $m,n,t,\qb$ tels que dans la d\'efinition \ref{defwm}.
Alors : 
$$\winfpr({\mathcal D}_{m,n,t}^{\qb}(k))=\wsuppr({\mathcal D}_{m,n,t}^{\qb}(k))=2m.$$
\end{propb}
\begin{demo}
Par d\'efinition on a w-supdeg$({\mathcal D}_{m,n,t}^{\qb}(k))\leq 2m$.
Mais d'apr\`es le th\'eor\`eme \ref{thwdeg} et la remarque pr\'ec\'edente, on a 
$$2m=\winf({\mathcal D}_{m,n,t}^{\qb}(k))\leq\wsup({\mathcal D}_{m,n,t}^{\qb}(k)).$$
\end{demo}


\medskip

A titre d'illustration, montrons sur l'exemple \'etudi\'e ci-dessus du corps de fonctions rationnelles $\Frac(S_{2,2}^{q}(k))$ comment un r\'esultat de ``d\'etressage par plongement'' permet le calcul du w-degr\'e sup\'erieur.

\subsubsection{Retour \`a l'exemple \ref{plg4}.}\label{wsup4} \label{eqwmsc4}

\begin{propb} 
Soient $q\in k^*$ non racine de l'unit\'e, et $\L=\left({\tiny \begin{array}{cc}
1& q\\
q^{-1}&1
\end{array}}\right)$.
Alors
$$\wsuppr(\Frac(S_{2,2}^{\L}(k)))=4.$$
\end{propb}
\begin{demo}
 La proposition \ref{plg4} montre que $\wsup(\Frac(S_{2,2}^{\L}(k)))\leq 4$.

Pour la r\'eciproque, donnons-nous un corps de Weyl mixte
 ${\mathcal D}_{m,n,t}^{\qb}(k)$  contenant   $\Frac(S_{2,2}^{\L}(k))$.
Montrons qu'alors  le corps de Weyl mixte ${\mathcal D}_{m,n+1,t}^{(\qb,q)}(k)$ contient le corps de Weyl ${\mathcal D}_2(k)=\Frac(A_2(k))$.
Pour cela, remarquons que
le corps ${\mathcal D}_{m,n+1,t}^{(\qb,q)}(k)$ est obtenu \`a partir de ${\mathcal D}_{m,n,t}^{\qb}(k)$ en  ajoutant au graphe de commutation donn\'e apr\`es la d\'efinition \ref{defwm} une ar\^ete quantique 
\begin{center}
\input{aq3.pstex_t}
\end{center}
Or  ${\mathcal D}_{m,n,t}^{\qb}(k)$ est suppos\'e contenir $\Frac(S_{2,2}^{\L}(k))$, engendr\'e par $x_1,x_2,y_1,y_2$ et les relations :
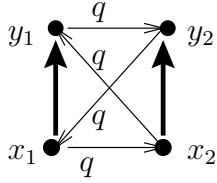
\begin{center}
\input{graf1.pstex_t}
\end{center}
Par ailleurs les g\'en\'erateurs $u_{n+1}$ et $v_{n+1}$ dans ${\mathcal D}_{m,n+1,t}^{(\qb,q)}(k)$ commutent \`a tous les \'el\'ements de ${\mathcal D}_{m,n,t}^{\qb}(k)$, et 
on v\'erifie  facilement que les \'el\'ements 
$$x'_1=v_{n+1}^{-1}x_1,\ y'_1=y_1v_{n+1},\ x'_2=u_{n+1}^{-1}x_2,\ y'_2=y_2u_{n+1}$$
 engendrent dans ${\mathcal D}_{m,n+1,t}^{(\qb,q)}(k)$ une image homomorphe (donc isomorphe) de l'alg\`ebre de Weyl $A_2(k)$.
 Donc  ${\mathcal D}_{m,n+1,t}^{(\qb,q)}(k)$ contient un sous-corps isomorphe au 
  corps de Weyl ${\mathcal D}_2(k)$.
Ceci prouve, d'apr\`es le th\'eor\`eme \ref{thwdeg}, que $m\geq 2$.
On en conclut par d\'efinition m\^eme du w-degr\'e sup\'erieur que :
$$\wsup(\Frac(S_{2,2}^{\L}(k)))\geq 4.$$
\end{demo}

\medskip

Comme on l'a annonc\'e au d\'ebut du paragraphe, ce calcul permet par exemple de retrouver le r\'esultat suivant, prouv\'e par une autre m\'ethode en \cite{AD3}.

\begin{corolb}
Soit $q\in k^*$ non racine de l'unit\'e, et $\L=\left({\tiny \begin{array}{cc}
1& q\\
q^{-1}&1
\end{array}}\right)$. Alors  le corps $\Frac(S_{2,2}^{\L}(k))$
n'est  $k$-isomorphe \`a aucun corps de Weyl mixte ${\mathcal D}_{m,n,t}^{\qb}(k)$.
\end{corolb}
\begin{demo}
Supposons que $\Frac(A_{2,2}^{\L}(k))$ soit  $k$-isomorphe \`a un corps de Weyl mixte ${\mathcal D}_{m,n,t}^{\qb}(k)$.
D'apr\`es la proposition \ref{heyhey},
l'\'egalit\'e des GK-degr\'es de transcendance impliquerait $4=2m+2n+t$.
Mais, d'apr\`es les propositions \ref{wwcwm} et  \ref{wsup4}, 
 l'\'egalit\'e des w-degr\'es sup\'erieurs impliquerait alors $4=2m$, donc $n=t=0$.
Le corps   ${\mathcal D}_{m,n,t}^{\qb}(k)$ serait alors le corps de Weyl ${\mathcal D}_2(k)$, ce qui est absurde puisque $G(\Frac(S_{2,2}^{\L}(k)))\neq \{1\}$.
\end{demo}

\subsection{Calcul du w-degr\'e sup\'erieur dans le cas g\'en\'eral.}

Nous allons  g\'en\'eraliser la m\'ethode de ``d\'etressage par plongement'' d\'ecrite sur l'exemple pr\'ec\'edent  afin de calculer le w-degr\'e sup\'erieur des corps de fonctions rationnelles mixtes crois\'es de dimension quelconque.

\begin{slemme} \label{ll1}
Soient $n\geq 1$ un entier,
et $\L\in M_n(k^*)$ une matrice multiplicativement antisym\'etrique. Alors il existe 
des entiers positifs $r,t$ v\'erifiant $2r+t=n(n-1)$, et un $r$-uplet $\qb\in (k\setminus\{0,1\})^r$  tels que  :
\begin{enumerate}
\item
l'espace affine quantique ${\mathcal O}_{\L}(k^n)$ se plonge dans l'alg\`ebre ${\mathcal A}_{0,r,t}^{\qb}(k)$ d\'efinie en \ref{defwm} ;
\item
et donc  $\Frac({\mathcal O}_{\L}(k^n))$  se plonge dans ${\mathcal D}^{\qb}_{0,r,t}(k)$.
\end{enumerate}
\end{slemme}
\begin{demo}
Consid\'erons l'alg\`ebre $A$  engendr\'ee sur $k$ par $n(n-1)$ g\'en\'erateurs  :

$$u^1_2,v^1_2, u^1_3,v^1_3,\ldots, u^1_n,v^1_n, u^2_3,v^2_3,\ldots, u^2_n,v^2_n,\ldots, u^{n-1}_n,v^{n-1}_n,$$
 soumis aux relations
\begin{equation} \label{relmix}
u^i_jv^i_j=\l_{i,j}v^i_ju^i_j, \textrm{ et } [u^i_j,v^k_t]= [u^i_j,u^k_t]= [v^i_j,v^k_t]=0 \textrm{ si } (i,j)\neq(k,t).
\end{equation}
 On v\'erifie ais\'ement que $A$ est une alg\`ebre ${\mathcal A}_{0,r,t}^{\qb}(k)$ au sens de la d\'efinition \ref{defwm},  dont le corps de fractions est un corps de Weyl mixte (en fait ici purement quantique) ${\mathcal D}^{\qb}_{0,r,t}(k)$, o\`u $t$ est le double du nombre de $\l_{i,j}$ \'egaux \`a 1, o\`u $r$ est le  nombre de $\l_{i,j}$ distincts de 1, et o\`u les $q_i$ dans $\qb=(q_1,\dots,q_r)$ sont pr\'ecis\'ement les $\l_{i,j}$ (avec $i<j$) distincts de 1.

On d\'efinit alors un morphisme de ${\mathcal O}_{\L}(k^n)=k_{\L}[y_1,\ldots,y_n]$ dans $A$ en posant :

$$
\begin{array}{l}
 \phi(y_i)=v^1_i\ldots v^{i-1}_iu^i_{i+1}\ldots u^i_n \ \textrm{ pour }\ 2\leq i\leq n-1,\\
  \phi(y_1)= u^1_{2}\ldots u^1_n,\ \textrm{ et }\ \phi(y_n)=v^1_n\ldots v^{n-1}_n.
\end{array}$$

Pour montrer  que $\phi$  est injectif on montre que pour tout 
$1\leq i\leq n$ le morphisme $\phi_i$ d\'efini comme la restriction de $\phi$ \`a la sous-alg\`ebre de $k_{\L}[y_1,\ldots,y_n]$ engendr\'ee par $y_1,\ldots,y_i$ est injectif. Notons $k_{\L_i}[y_1,\ldots,y_i]$ cette sous-alg\`ebre.
Pour $i=1$, 
 l'injectivit\'e de $\phi_1$ est claire.
Soit $i\geq 2$, et
supposons $\phi_1,\ldots,\phi_{i-1}$ injectifs. Notons $P=\ssum_{a=0}^d P_a(y_1,\ldots,y_{i-1})y_i^a$ un \'el\'ement non nul de $k_{\L_i}[y_1,\ldots,y_i]$, avec $P_a(y_1,\ldots,y_{i-1})\in k_{\L_{i-1}}[ y_1,\ldots,y_{i-1}]$, et $P_d\neq 0$. 
Alors 
$$\phi_i(P) =\sum_{a=0}^d \phi_i(P_a(y_1,\ldots,y_{i-1})) (v^1_i)^a \ldots (v^{i-1}_i)^a(u^i_{i+1})^a\ldots (u^i_n)^a.$$
Puisque $v^{i-1}_i$ n'appara\^\i t dans aucun des $\phi_i(y_1),\ldots,\phi_i(y_{i-1})$, le terme de plus haut degr\'e en $v^{i-1}_i$ est  $\phi_i(P_d(y_1,\ldots,y_{i-1}))(v^1_i)^d \ldots (v^{i-1}_i)^d(u^i_{i+1})^d\ldots (u^i_n)^d$.
Or par hypoth\`ese de r\'ecurrence, $\phi_i(P_d(y_1,\ldots,y_{i-1}))=\phi_{i-1}(P_d(y_1,\ldots,y_{i-1}))\neq 0$, donc $\phi_i(P)\neq 0$, ce qui montre l'injectivit\'e de $\phi_i$.
Ainsi par it\'erations $\phi_n=\phi$ est injectif. Ceci montre le point (1).
Pour le point (2),
on \'etend  $\phi$ en un plongement de $k_{\L}(y_1,\ldots,y_n)$ dans $\Frac(A)={\mathcal D}^{\qb}_{0,r,t}(k)$.
\end{demo}

\medskip

Illustrons ce lemme en dimension 3.
Soit $\L\in M_3(k^*)$ multiplicativement antisym\'etrique.
L'espace affine quantique $k_{\L}[y_1,y_2,y_3]$ est engendr\'e par $y_1,y_2,y_3$ avec les relations :

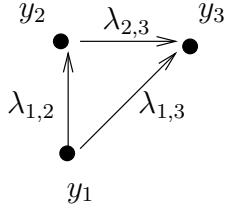
\begin{center}
\input{ptq2.pstex_t}
\end{center}

L'alg\`ebre $A$ de la preuve est engendr\'ee par $u_2^1,v_2^1,u_3^1,v_3^1,u_3^2,v_3^2$ avec les relations :

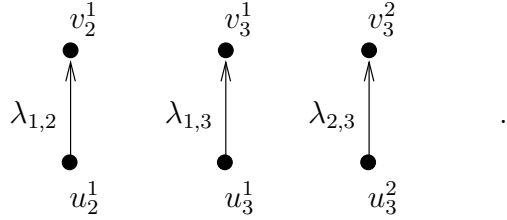
\begin{center}
\input{plg3.pstex_t}
\end{center}

On d\'efinit alors un morphisme injectif de $k_{\L}[y_1,y_2,y_3]$ dans $A$ par
$y_1\mapsto u^1_2u^1_3$, $y_2\mapsto v^1_2u^2_3$, $y_3\mapsto v^1_3v^2_3$.

\medskip

Montrons maintenant comment ``d\'etresser'' un graphe  comprenant \`a la fois des ar\^etes eul\'eriennes et des ar\^etes quantiques.

\begin{slemme} \label{ll2}
Soit $q\in k\setminus\{0,1\}$. Alors le corps mixte crois\'e  $\Frac(S_{2,1}^q(k))$ se plonge dans le corps de Weyl  mixte ${\mathcal D}_{1,1,0}^q(k)$.
\end{slemme}
\begin{demo}
L'alg\`ebre $T^q_{2,1}(k)$ d\'efinie en \ref{arres} est engendr\'ee sur $k$ par $y_1,y_2,w_1$ avec les relations :
\begin{center}
$y_1y_2=qy_2y_1$, $y_1w_1=(w_1-1)y_1$ et $y_2w_1=w_1y_2$,
\end{center}
Par ailleurs, consid\'erons l'alg\`ebre $A=U(k)\otimes {\mathcal O}_q(k^2)$ engendr\'e sur $k$ par $w,y,u,v$, avec les relations :

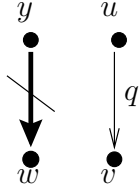
\begin{center}
\input{graf4.pstex_t}
\end{center}

 On d\'efinit un morphisme de $k$-alg\`ebres de $T^q_{2,1}(k)$ dans $A$
en envoyant $y_1$ sur $yu$, $y_2$ sur $v$, $w_1$ sur $w$. Les mon\^omes en  les g\'en\'erateurs formant une base de $k$-espace vectoriel de $A$, on montre comme en \ref{ll1} que ce morphisme est injectif.
 On \'etend alors ce morphisme en un plongement de ${\Frac} (T^Q_{2,1}(k))$ 
 dans le corps de Weyl mixte ${\mathcal D}_{1,1,0}^q(k)=\Frac(A)$.
\end{demo}

\medskip

Enon\c cons enfin le r\'esultat g\'en\'eral.

\begin{spropo}  \label{plgtdlambda}
Soient $n\geq 1$ et $0\leq r\leq n$ deux entiers, 
et $\L\in M_n(k^*)$ une matrice multiplicativement an\-ti\-sy\-m\'e\-tri\-que. Alors il existe des entiers positifs $s,t$ v\'erifiant  $n(n-1)\leq 2s+t\leq n(n-1)+r$, et un $s$-uplet $\qb\in(k\setminus\{0,1\})^s$ tels que le corps de fonctions rationnelles mixte crois\'e
 $\Frac(S^{\L}_{n,r}(k))$ se plonge dans le corps de Weyl mixte ${\mathcal D}_{r,s,t}^{\qb}(k)$.
\end{spropo}
\begin{demo}
Si $n=1$ l'alg\`ebre $S_{1,r}^{(1)}(k)$ vaut $k[y_1]$ ou $A_1(k)$ suivant que $r=0$ ou $r=1$, et il n'y a rien \`a d\'emontrer. Supposons donc $n\geq 2$.
Soient $y_1,\ldots,
y_n,w_1,\ldots,w_r$ les g\'en\'erateurs de la sous-alg\`ebre $T^{\L}_{n,r}(k)$ de $S_{n,r}^{\L}(k)$ d\'efinie en \ref{arres}.
On va montrer que  
le plongement consid\'er\'e au lemme \ref{ll1} pour la sous-alg\`ebre engendr\'ee par $y_1,\ldots,y_n$ s'\'etend en un plongement de $T^{\L}_{n,r}(k)$ tout entier dans une alg\`ebre $A$ dont le corps de fractions est un corps de Weyl mixte.
Soit $A$
l'alg\`ebre  engendr\'ee sur $k$ par $n(n-1)+r$ g\'en\'erateurs 
$$u^1_2,v^1_2, u^1_3,v^1_3,\ldots, u^1_n,v^1_n, u^2_3,v^2_3,\ldots, u^2_n,v^2_n,\ldots, u^{n-1}_n,v^{n-1}_n,t_1,\ldots,t_r,$$
avec les relations suivantes :
\begin{enumerate}
\item
  les $n(n-1)$ premiers g\'en\'erateurs $u^i_j$ et $v^k_l$ v\'erifient les relations (\ref{relmix})  d\'efinies en  \ref{ll1} ;
\item
pour $i<n$, l'\'el\'ement $t_i$ commute \`a tous les g\'en\'erateurs except\'e $u_{i+1}^i$, avec lequel il v\'erifie la relation $[t_i,u_{i+1}^i]=u_{i+1}^i$ ;
\item
enfin si $r=n$, alors $t_n$ commute \`a tous les autres g\'en\'erateurs except\'es $v^1_n$, avec lequel il v\'erifie $[t_n,v_{n}^1]=v_{n}^1$. 
\end{enumerate}

L'alg\`ebre ainsi d\'efinie est une extension it\'er\'ee de Ore en ces $n(n-1)+r$ g\'en\'erateurs,
 et admet donc un corps de fractions $K$.
On peut  d\'efinir un morphisme  de $T^{\L}_{n,r}(k)$ dans $A$, en compl\'etant le morphisme d\'efini au lemme \ref{ll1} par $w_i\mapsto t_i$. 
On d\'emontre comme en \ref{ll1} que ce morphisme est injectif.

Le graphe $\G$ associ\'e aux relations de commutation entre les g\'en\'erateurs de $A$ est constitu\'e :
\begin{itemize}
\item  de sommets d'o\`u ne part ni n'arrive aucune ar\^ete,  correspondant \`a des coefficients $\l_{i,j}$ valant 1,
\item de sommets reli\'es 2 \`a 2 par une ar\^ete quantique et d'o\`u ne part aucune autre ar\^ete que celle-ci, 
\item et de triangles du type :

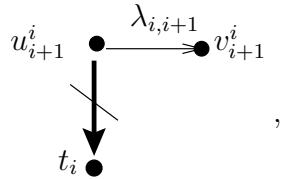
\begin{center}
\input{graf5.pstex_t}
\end{center}
dont les sommets ne sont reli\'es \`a aucun autre sommet de $\G$.
\end{itemize}

Ainsi $A$ est produit tensoriel d'une alg\`ebre de polyn\^omes commutatifs, d'un nom\-bre fini  de plans quantiques, et d'un nombre fini d'alg\`ebres $T_{2,1}^{\l_{i,i+1}}(k)$ d\'efinies par un graphe du type ci-dessus.

Si $\l_{i,i+1}=1$, alors $T_{2,1}^{\l_{i,i+1}}(k)$ est le produit tensoriel de $k[v^i_{i+1}]$ par une alg\`ebre $U(k)$ dont le corps de fractions est le corps  de Weyl $D_1(k)$. 
Si $\l_{i,i+1}\neq 1$, de fa\c con analogue \`a ce qui a \'et\'e fait au lemme \ref{ll2}, on plonge $T_{2,1}^{\l_{i,i+1}}(k)$ dans une
alg\`ebre produit tensoriel d'une alg\`ebre $U(k)$ par un plan quantique.

On plonge  ainsi  $A$ dans une alg\`ebre $B$ produit tensoriel de plans quantiques, d'alg\`ebres $U(k)$, et d'une alg\`ebre de polyn\^omes commutatifs.
On en d\'eduit un morphisme injectif de $T_{n,r}^{\L}(k)$ dans $B$.
On \'etend alors ce morphisme en un plongement du corps de fractions de $T_{n,r}^{\L}(k)$
 dans le corps de fractions de $B$, qui est par construction 
un corps de Weyl mixte.
\end{demo}

\medskip

Le lemme suivant est la derni\`ere \'etape vers le calcul du
 w-degr\'e sup\'erieur d'un corps de fonctions rationnelles mixte crois\'e.

\begin{slemme} \label{yeman}
Soient $n\geq 1$ et $0\leq r\leq n$ deux entiers, et
 $\L\in M_n(k^*)$ une matrice multiplicativement antisym\'etrique. Supposons qu'il existe des entiers $m,s,t$, et $\qb\in(k\setminus\{0,1\})^s$ tels que le corps de Weyl mixte ${\mathcal D}_{m,s,t}^{\qb}(k)$ admette un sous-corps $k$-isomorphe au corps de fonctions rationnelles mixte crois\'e
 $\Frac(S^{\L}_{n,r}(k))$. 
Alors $r\leq m$.
\end{slemme}
\begin{demo}
Soit ${}^t\L$ la matrice multiplicativement antisym\'etrique transpos\'ee de $\L$.
D'apr\`es le lemme \ref{ll1}, il existe une alg\`ebre 
 ${\mathcal A}_{0,s_1,t_1}^{\qb'}(k)$ au sens de \ref{defwm}  telle que ${\mathcal O}_{^t\L}(k^n)$ se plonge dans ${\mathcal A}_{0,s_1,t_1}^{\qb'}(k)$.
Pour  les donn\'ees $m,s,t,\qb$ de l'\'enonc\'e, consid\'erons l'alg\`ebre en $2m+2s+t+2s_1+t_1$ g\'en\'erateurs
  $$B={\mathcal A}_{m,s,t}^{\qb}(k)\otimes {\mathcal A}_{0,s_1,t_1}^{\qb'}(k)\cong {\mathcal A}_{m,s+s_1,t+t_1}^{(\qb,\qb')}(k),$$
 o\`u $(\qb,\qb')=(q_1,\ldots,q_s,q'_1,\ldots,q'_{s_1})\in (k\setminus\{0,1\})^{s+s_1}$.

Par hypoth\`ese $\Frac(S_{n,r}^{\L}(k))$ se plonge dans ${\mathcal D}_{m,s,t}^{\qb}(k)$ donc dans $\Frac(B)$.
En particulier, il existe dans $\Frac(B)$ des \'el\'ements $y_1,\ldots,y_n,w_1,\ldots,w_r$ engendrant sur $k$ une sous-alg\`ebre $k$-isomorphe \`a $T^{\L}_{n,r}(k)$, et des \'el\'ements $y'_1,\ldots,y'_n$ engendrant sur $k$ une sous-alg\`ebre $k$-isomorphe \`a ${\mathcal O}_{^t\L}(k^n)$, et tels que les $y'_i$ commutent aux $y_1,\ldots,y_n,w_1,\ldots,w_r$.
Posons  $Y_k=y_ky'_k$, et $X_k =Y_k^{-1} w_k$ dans $\Frac(B)$  pour tout $k\leq r$.
Les \'el\'ements $X_1,\ldots,X_r,Y_1,\ldots,Y_r$ engendrent dans $\Frac(B)$ une image homomorphe (donc isomorphe) de l'alg\`ebre de Weyl  $A_r(k)$.
On conclut que   ${\mathcal D}_r(k)$ se plonge dans $\Frac(B)={\mathcal D}_{m,s+s_1,t+t_1}^{(\qb,\qb')}(k)$.
Il d\'ecoule alors du th\'eor\`eme  \ref{thwdeg} que $r\leq m$.
\end{demo}

\begin{stheoreme} \label{gdwdeg}
Soit $\L\in M_n(k^*)$ une matrice multiplicativement an\-ti\-sy\-m\'e\-tri\-que, et $r\leq n$.
Alors 
$$\wsuppr(\Frac(S^{\L}_{n,r}(k)))=2r.$$
\end{stheoreme}
\begin{demo}
D'apr\`es  le  lemme  \ref{yeman} on a  $\wsup(\Frac(S^{\L}_{n,r}(k)))\geq 2r$, et la proposition \ref{plgtdlambda} montre l'in\'egalit\'e inverse.
\end{demo}

\medskip

{\bf Remarque.}
On a ainsi caract\'eris\'e $r$ comme un invariant rationnel
 des alg\`ebres
 $S^{\L}_{n,r}(k)$, ind\'ependant de la pr\'esentation choisie.
Rappelons que, si l'on voit $S_{n,r}^{\L}(k)$ comme une alg\`ebre $\AG$ d\'efinie \`a partir 
  d'un graphe $qW$
 admissible $\G$ au sens de la section \ref{secgraf},
 alors $r$ n'est autre que 
 la moiti\'e du rang de la matrice $P(\G)$ codant les ar\^etes de Weyl dans $\G$.
Pour certaines classes d'alg\`ebres rationnellement \'equivalentes \`a des alg\`ebres $S_{n,r}^{\L}(k)$, il est int\'eressant de pouvoir calculer $r$ directement en fonction des param\`etres d\'efinissant ces alg\`ebres.
C'est le cas par exemple des alg\`ebres de Weyl quantiques multiparam\'etr\'ees, comme on le verra \`a la section \ref{secawq}.

\section{Equivalence rationnelle des alg\`ebres  mixtes crois\'ees.}

Commen\c cons par donner une condition n\'ecessaire pour l'\'equivalence rationnelle de deux alg\`ebres polynomiales mixtes crois\'ees.

\begin{theoreme} \label{bientot}
Soient $n,n'\geq 1$ et $r\leq n$, $r'\leq n'$ des entiers.
Soient deux matrices $\L\in M_n(k^*)$, et $\L'\in M_{n'}(k^*)$ multiplicativement antisym\'etriques.
 Supposons que les  deux corps   $\Dr$ et $\Frac(S^{\L'}_{n',r'}(k))$ soient $k$-isomorphes. Alors :
\begin{enumerate}
\item
$n=n'$ et $r=r'$.
\item
Si de plus le tore quantique  ${\mathcal O}_{\L}({k^*}^n)$ est simple, alors les tores quantiques
${\mathcal O}_{\L}({k^*}^n)$ et ${\mathcal O}_{\L'}({k^*}^{n'})$ sont $k$-isomorphes.
\end{enumerate}
\end{theoreme}
\begin{demo}
Par \'egalit\'e des GK-degr\'es de transcendance  on a d'apr\`es la proposition \ref{gkdsnrl} l'\'egalit\'e $n+r=n'+r'$. Par ailleurs le th\'eor\`eme  \ref{gdwdeg} prouve  que $r=r'$, donc $n=n'$. Ainsi $\L$ et $\L'$ ont la m\^eme taille, et  le point 2   r\'esulte alors du  lemme \ref{kifoqr} et de la proposition \ref{happy}.
\end{demo}

\medskip

La question se pose naturellement de la r\'eciproque du point 2 du th\'eor\`eme \ref{bientot}.
Il est assez facile de v\'erifier qu'une telle r\'eciproque est vraie dans le contexte des corps de Weyl mixtes (voir le th\'eor\`eme 3.1.3.3 de \cite{these}).
Nous n'avons de r\'eponse compl\`ete pour les corps de fonctions rationnelles mixtes crois\'es $\Frac(S_{n,r}^{\L}(k))$ que dans le cas semi-classique ($n=r$), que nous allons maintenant d\'evelopper.
Observons d'abord que les corps de fonctions rationnelles mixtes crois\'es semi-classiques $\Frac(S_{n,n}^{\L}(k))$ ne sont des corps de Weyl mixtes que si ce sont des corps de Weyl classiques.
C'est l'objet de la proposition suivante, qui g\'en\'eralise le corollaire \ref{eqwmsc4}.

\begin{propo}
Soient $n\geq 1$ un entier, et  $\L\in M_{n}(k^*)$ multiplicativement antisym\'etrique. 
Supposons qu'il existe des entiers positifs $m',n',t'$ et un $n'$-uplet $\qb\in(k\setminus\{0,1\})^{n'}$ tels que  $\Frac(S^{\L}_{n,n}(k))$ soit $k$-isomorphe au corps de Weyl mixte ${\mathcal D}_{m',n',t'}^{\qb}(k)$. Alors $n=m'$, et $n'=t'=0$, c'est-\`a-dire que $\Frac(S^{\L}_{n,n}(k))$ est le  corps de Weyl classique ${\mathcal D}_n(k)=\Frac(A_n(k))$.
\end{propo}
\begin{demo}
Par \'egalit\'e des GK-degr\'es de transcendance  on a : $2n=2m'+2n'+t'$. 
Par ailleurs d'apr\`es la proposition \ref{wwcwm} et le th\'eor\`eme \ref{gdwdeg}, l'\'egalit\'e du w-degr\'e sup\'erieur des deux corps implique $2n=2m'$, ce qui prouve la proposition.
\end{demo}

\medskip

Par \'egalit\'e des GK-degr\'es de transcendance, 
la r\'eciproque au th\'eor\`eme \ref{bientot} dans le cas semi-classique n'a de sens que pour $n=n'$. Elle s'\'enonce alors comme suit.

\begin{theoreme}[Cas semi-classique] \label{pigne}
Soient $n \in {\mathbb N}$, et
 $\L,\L' \in M_n(k^*)$  multiplicativement antisym\'etriques.
Supposons que le  tore quantique ${\mathcal O}_{\L}({k^*}^n)$ est simple. Alors  les corps $\Frac(S^{\L}_{n,n}(k))$ et $\Frac(S^{\L'}_{n,n}(k))$ sont $k$-isomorphes si et seulement si les tores ${\mathcal O}_{\L}({k^*}^n)$ et ${\mathcal O}_{\L'}({k^*}^{n})$ sont  $k$-i\-so\-mor\-phes.
\end{theoreme}
\begin{demo}
La condition est n\'ecessaire par le th\'eor\`eme \ref{bientot}.

R\'eciproquement, soit $\Phi$ un isomorphisme   du tore quantique  $k_{\L}[y_1^{\pm 1},\ldots,y_n^{\pm 1}]$ sur
$k_{\L'}[{y'}_1^{\pm 1},\ldots,{y'}_n^{\pm 1}]$. Il existe d'apr\`es la proposition \ref{MP} une matrice $H=(h_{i,j})_{1\leq i,j\leq n}$ dans $GL_n({\mathbb Z})$ telle que pour tout $i$, $\Phi(y_i)=\a_i{y'_1}^{h_{1,i}}\ldots {y'_n}^{h_{n,i}}$, avec $\a_i \in k^*$. Notons $H^{-1} = (h'_{i,j})_{1\leq i,j\leq n}$.
Notons $w_1,\ldots,w_n,y_1,\ldots,y_n$ les g\'en\'erateurs sur $k$ de la sous-alg\`ebre $T_{n,n}^{\L}(k)$ de $S_{n,n}^{\L}(k)$ d\'efinie en \ref{arres}.1, et $w'_1,\ldots,w'_n,y'_1,\ldots,y'_n$ les g\'en\'erateurs sur $k$ de $T_{n,n}^{\L'}(k)$.
On v\'erifie alors ais\'ement qu'il existe un isomorphisme $\widehat \Phi$ de $\what T^{\L}_{n,n}(k)$ sur $\what T^{\L'}_{n,n}(k)$, d\'efini pour tous $1\leq i,j\leq n$  par $\widehat \Phi(y_i)={y'_1}^{h_{1,i}}\ldots {y'}_n^{h_{n,i}}$ et $\widehat \Phi(w_j) = h'_{j,1}w'_1 + \ldots +h'_{j,n}w'_n$.
Cet isomorphisme s'\'etend aux corps de fractions, ce qui termine la preuve.
\end{demo}

\medskip

{\bf Remarque.}
 Le th\'eor\`eme  \ref{pigne} permet de retrouver certains r\'esultats de \cite{Retq}. En effet, sa preuve  montre en fait, pour les alg\`ebres d'op\'erateurs diff\'erentiels eul\'eriens sur les tores quantiques simples (c'est-\`a-dire les alg\`ebres $\what T_{n,n}^{\L}(k)$ avec les notations de \ref{arres}.1), que les conditions suivantes sont \'equivalentes :
\begin{enumerate}
\item $\what T_{n,n}^{\L}(k)$ est $k$-isomorphe \`a $\what T_{n,n}^{\L'}(k)$;
\item $\what T_{n,n}^{\L}(k)$ est rationnellement \'equivalente  \`a $\what T_{n,n}^{\L'}(k)$;
\item les tores quantiques simples ${\mathcal O}_{\L}({k^*}^n)$ et ${\mathcal O}_{\L'}({k^*}^n)$  sont $k$-isomorphes;
\item les tores quantiques simples ${\mathcal O}_{\L}({k^*}^n)$ et ${\mathcal O}_{\L'}({k^*}^n)$  sont rationnellement \'equivalents.
\end{enumerate}

Notons enfin que 
dans le cas uniparam\'etr\'e, m\^eme \emph{sans l'hypoth\`ese de simplicit\'e du tore quantique sous-jacent},
on peut d\'emontrer un r\'esultat similaire, gr\^ace \`a la notion suivante de matrice canonique, utilis\'ee  par A.N. Panov dans \cite{Pa1}.

\begin{defin} \label{defcan}
Soient deux entiers $n\geq 1$  et $0 \leq s\leq n/2$.
Soit $(c_1,\ldots,c_s)$ une famille d'entiers strictement positifs tels que $c_i$ divise $c_{i+1}$ pour tout $i$.
On appelle matrice canonique antisym\'etrique de taille $n$ associ\'ee \`a $(c_1,\ldots,c_s)$,
la matrice diagonale par blocs 
$$C^n(c_1,\ldots,c_s)=\textrm{Diag}[C_1,\ldots,C_s,0,\ldots,0],\ \textrm{ o\`u }\  C_i={\tiny \left( \begin{array}{cc}  
0      & c_i   \cr
-c_i   & 0     \cr
\end{array}
\right)}.$$
\end{defin}

\begin{theoreme}[Cas uniparam\'etr\'e]
Soient $n\geq 1$ un entier, et $q\in k^*$ non racine de l'unit\'e.
\begin{itemize}
\item[(i)] Soit $S=(s_{i,j})\in M_n({\mathbb Z})$ une matrice antisym\'etrique de rang $2s$, et posons $\L=\left(q^{s_{i,j}}\right)_{i,j} \in M_n(k^*)$. Alors il existe
une matrice antisym\'etrique canonique $C=(c_{i,j})=C^n(d_1,\ldots,d_s)$ \'equivalente \`a $S$, telle que 
  $\Frac(S^{\L}_{n,n}(k))$ est $k$-isomorphe \`a $\Frac(S^{\L'}_{n,n}(k))$, o\`u $\L'=(q^{c_{i,j}})_{i,j}$.
\item[(ii)] Soient $s,s'\leq n/2$ deux entiers, et deux matrices antisym\'etriques canoniques  $C=C^n(d_1,\ldots,d_s)$ et $C'=C^n(d'_1,\ldots,d'_{s'})$.
 Alors les corps gauches $\Frac(S^{\L}_{n,n}(k))$ et $\Frac(S^{\L'}_{n',n'}(k))$ sont $k$-isomorphes si et seulement si $n=n'$, et $C=C'$.
\end{itemize}
\end{theoreme}
\begin{demo}
On adapte les preuves de \cite{Pa1}, en s'assurant qu'on peut remplacer le corps de base $k$ par le corps de Weyl $D_n(k)$.
\end{demo}

\section{Applications aux alg\`ebres de Weyl quantiques multiparam\'etr\'ees.}
\label{secawq}

Un exemple particuli\`erement significatif d'alg\`ebres rationnellement \'equivalentes \`a des alg\`ebres polynomiales mixtes crois\'ees est celui des alg\`ebres de Weyl quantiques multiparam\'etr\'ees.
Notons que dans le cas ``purement quantique'', cette \'equivalence rationnelle a \'et\'e \'etablie dans \cite{AD1}.
Il existe dans la litt\'erature de tr\`es nombreuses fa\c cons d'introduire les alg\`ebres de Weyl quantiques associ\'ees \`a  la d\'efinition d'un calcul diff\'erentiel sur diverses d\'eformations de l'espace affine (\cite{WZ}, \cite{Malt}, \cite{De}, \cite{GZ}...). 
Toutes aboutissent \`a des formes particuli\`eres  de la d\'efinition g\'en\'erale que nous donnons ci-dessous (d'apr\`es \cite{Malt}; voir aussi \cite{BG}, d\'efinition I.2.6, et la section 1.3 de \cite{these}).
Ces alg\`ebres ont fait l'objet de nombreuses \'etudes du point de vue de la th\'eorie des anneaux (voir par exemple   les articles  \cite{AD1}, \cite{Jslqwa}, \cite{GZ},
 \cite{FKK},
 \cite{F} et leurs bibliographies, ainsi que le chapitre 5 de \cite{these}).

\begin{defin} Soit $n \geq 1$ un entier, soit $\L \in M_n(k^*)$ multiplicativement antisym\'etrique,
et soit $ \qb = (q_1,\ldots,q_n) \in (k^*)^n$. L'alg\`ebre de Weyl
multiparam\'etr\'ee $\An(k)$  est la $k$-alg\`ebre engendr\'ee   par $2n$ g\'en\'erateurs $x_1,y_1,\ldots,x_n,y_n$ soumis aux relations ci-dessous, pour $1\leq i<j\leq n$   :
\begin{equation} \label{relan}
\begin{array}{c}
y_i y_j = \l_{i,j} y_j y_i,\\
x_i x_j = q_i \l_{i,j} x_j x_i,\\
x_i y_j = \l_{i,j}^{-1} y_j x_i,\\
x_j y_i = q_i \l_{i,j} y_i x_j,\\
x_j y_j = 1 + \sum_{1\leq k<j} (q_k - 1)y_k x_k + q_j y_j x_j.
\end{array}
\end{equation}
\end{defin}

Comme pour les alg\`ebres polynomiales mixtes crois\'ees on peut v\'erifier  que la GK-dimension de $\An(k)$ et le GK-degr\'e de transcendance de $\Frac(\An(k))$ valent $2n$.

Les alg\`ebres de Weyl quantiques multiparam\'etr\'ees d\'efinies ci-dessus, ou leur version ``alternative" ${\mathcal A}_n^{\qb,\L}(k)$  pr\'esent\'ee dans \cite{AJ}, admettent toujours  une localisation commune avec des alg\`ebres $\Alr$.
C'est l'objet de la proposition suivante.

\begin{propo} \label{eqansn}
Soient $n\geq 1$ un entier, $\qb=(q_1,\ldots,q_n)\in(k^*)^n$ et $\L\in M_n(k^*)$ multiplicativement antisym\'etrique.
Soient $r$ le nombre de $q_i$ \'egaux \`a 1, et $q_{i_1},\ldots$, $q_{i_{n-r}}$  les $q_i$ diff\'erents de 1. Soit  $\L'\in M_{2n-r}(k^*)$ multiplicativement antisym\'etrique
$$\L'=\left(\begin{array}{cc}
(1) & M\\
 M'&\L
\end{array}\right),\ \textrm{ avec} $$

\begin{tabular}{l}
$M=(m_{k,t})\in M_{n-r,n}(k^*)$ d\'efinie par $m_{k,t}=q_{i_k}$ si $t=i_k$ et $m_{k,t}=1$ sinon,\\
$M'=(m'_{k,t})\in M_{n,n-r}(k^*)$ d\'efinie par $m'_{k,t}=q_{i_t}^{-1}$ si $k=i_t$ et $m'_{k,t}=1$ sinon.
\end{tabular}

Alors $\Frac(\An(k))$ est $k$-isomorphe \`a $\Frac(S_{2n-r,r}^{\L'}(k))$.
\end{propo}
\begin{demo}
Dans $\An(k)$ soit $z_i=x_iy_i-y_ix_i$ pour tout $1\leq i\leq n$ (voir \cite{Jslqwa}, et la section 1.3 de \cite{these}). Ce sont des \'el\'ements normaux, et on note $\Bn(k)$ le localis\'e de $\An(k)$ en la partie multiplicative engendr\'ee par les $z_i$.
Consid\'erons  les \'el\'ements $x'_j=z_{j-1}^{-1}x_j$   et $z'_i=z_{i-1}^{-1}z_i$ de $\Bn(k)$.
Alors dans $\Bn(k)$ les \'el\'ements 
$$y_{j_1},\ldots,y_{j_r},y_{i_1},\ldots, y_{i_{n-r}},x'_{j_1},\ldots,x'_{j_r},z'_{i_1},z'_{i_{n-r}},$$
o\`u $j_1<\ldots<j_r$ sont les indices tels que $q_{j_k}=1$ pour $1\leq k\leq r$, et $i_1<\ldots<i_{n-r}$ 
sont les indices tels que $q_{i_k}\neq1$ pour $1\leq k\leq n-r$,
 v\'erifient les relations :
\begin{equation} \label{marilyn}
\begin{array}{cl}
y_iy_j=\l_{i,j}y_jy_i, & \textrm{pour tous } \ 1\leq i,j \leq n,\\
y_ix'_{j_k}=\l_{i,j_k}^{-1}x'_{j_k}y_i, &\textrm{pour tous } \ 1\leq i \leq n,\ 1\leq k\leq r, i\neq j_k\\
x'_{j_k}y_{j_k}=y_{j_k}x'_{j_k}+1, &\textrm{pour tous } \ 1\leq k\leq r\\
y_jz'_{i_k}=z'_{i_k}y_j, &\textrm{pour tous } \ 1\leq j \leq n,\ 1\leq k\leq n-r,\ j\neq i_k\\
y_{i_k}z'_{i_k}=q_{i_k}z'_{i_k}y_{i_k}, &\textrm{pour tous } \ 1\leq k\leq n-r,\\
x'_{j_k}z'_{i_l}=z'_{i_l}x'_{j_k}, &\textrm{pour tous }\ 1\leq k\leq r, \ 1\leq l\leq n-r.
\end{array}
\end{equation}
On v\'erifie ais\'ement que les mon\^omes en ces \'el\'ements sont
lin\'eairement in\-d\'e\-pen\-dants.
Ceux-ci engendrent  donc une sous-alg\`ebre de $\Bn(k)$ $k$-isomorphe \`a l'alg\`ebre polynomiale mixte crois\'ee $\smash{S_{2n-r,r}^{\L'}(k)}$.
Par ailleurs  leur d\'efinition  implique que le corps de fractions de cette sous-alg\`ebre  est $\Frac(\An(k))$ tout entier.
\end{demo}

\begin{corol}
Soit $\An(k)$ une alg\`ebre de Weyl quantique multiparam\'etr\'ee.
On suppose que, pour tout $1\leq i\leq n$, on a $q_i$=1 ou $q_i$ non racine de l'unit\'e.
Alors ${\mathcal Z}(\Frac(\An(k)))=k.$
\end{corol}
\begin{demo}
On a vu \`a la proposition \ref{eqansn}
 que le corps de fractions de $\An(k)$ est  le corps de fractions de l'alg\`ebre polynomiale mixte crois\'ee engendr\'ee par les g\'en\'erateurs 
$$y_{j_1},\ldots,y_{j_r},y_{i_1},\ldots, y_{i_{n-r}},x'_{j_1},\ldots,x'_{j_r},z'_{i_1},z'_{i_{n-r}},$$
avec les relations (\ref{marilyn}).
Avec les notations de \ref{eqansn}, la proposition \ref{centre} montre que ${\mathcal Z}(\Frac(\An(k)))$ est  l'intersection de ${\mathcal Z}(k_{\L'}(y_1,\ldots,y_n,z'_{i_1},\ldots,z'_{i_{n-r}}))$ avec le sous-corps de $\Frac(\An(k))$ engendr\'e par $y_{i_1},\ldots,y_{i_{n-r}},z'_{i_1},\ldots,z'_{i_{n-r}}$.
Puisqu'aucun des $q_{i_1},\ldots,q_{i_{n-r}}$ n'est racine de l'unit\'e, on montre facilement que le centre du tore quantique  ${\mathcal O}_{\L'}({k^*}^{2n-r})$ ne comporte aucun mon\^ome en les g\'en\'erateurs $y_{i_1},\ldots,y_{i_{n-r}},z'_{i_1},\ldots,z'_{i_{n-r}}$. On en d\'eduit avec la proposition 2.8 de \cite{Ca1}  qu'il en est de m\^eme du centre du corps  $k_{\L'}(y_1,\ldots,y_n,z'_{i_1},\ldots,z'_{i_{n-r}})$, et donc que  ${\mathcal Z}(\Frac(\An(k)))=k$.
\end{demo}

\begin{propo} \label{wdegdnql}
Soient $n\geq 1$ un entier,
 $\L\in M_n(k^*)$ une matrice multiplicativement antisym\'etrique, et $\qb=(q_1,\ldots,q_n)\in(k^*)^n$.
Notons $r$ le nombre de $q_i$ \'egaux \`a 1.
Alors 
$$\wsuppr(\Frac(\An(k)))=2r.$$
\end{propo}
\begin{demo}
Ceci d\'ecoule  du th\'eor\`eme \ref{gdwdeg} et de la proposition \ref{eqansn}.
\end{demo}

On obtient enfin les r\'esultats suivants concernant les corps de fractions des alg\`ebres de Weyl quantiques multiparam\'etr\'ees.

\begin{theoreme} Soient $n,n'\geq 1$ deux entiers,
 $\qb\in(k^*)^n$, $\qb'\in(k^*)^{n'}$,
 $\L\in M_n(k^*)$ et $\L'\in M_{n'}(k^*)$ multiplicativement antisym\'etriques.\\
 Si  ${\rm Frac}(\An(k))\cong {\rm Frac}(A^{\qb',\L'}_{n'}(k))$, alors :
\begin{itemize}
\item
$n=n'$ et
 $\qb$ et $\qb'$ ont le m\^eme nombre $r$ d'\'el\'ements \'egaux \`a 1.\\[.2cm]
\item
Si de plus le tore quantique  ${\mathcal O}_{\widetilde\L}({k^*}^{2n-r})$ est simple, avec les notations de la proposition \ref{eqansn}, alors: 
$${\mathcal O}_{\widetilde\L}({k^*}^{2n-r})\cong {\mathcal O}_{\widetilde{\L'}}({k^*}^{2n-r}).$$
\end{itemize}
\end{theoreme}
\begin{demo}
On applique le th\'eor\`eme \ref{bientot}, {\it via} la proposition \ref{eqansn}.
\end{demo}

\medskip

{\bf Remarque.} Si tous les $q_i$ sont \'egaux \`a 1, c'est-\`a-dire que $r=n$ alors $A_n^{\bar1,\L}(k)=S_{n,n}^{\L}(k)$, et le th\'eor\`eme \ref{pigne}  s'applique directement.


\begin{thebibliography}{99}
                        
\bibitem{AJ}
 M. Akhavizadegan, D.A. Jordan,
Prime ideals of quantized Weyl algebras,
 Glasgow Math. J.
{ 38}
(1996), no. 3,
283--297.



\bibitem{AD1}
J. Alev, F. Dumas, 
Sur le corps de fractions de certaines alg\`ebres quantiques,
 J. Algebra
{ 170} (1994),
no. 1,
 229--265.
 

\bibitem{AD3}
$\soul\esp$,
Corps de Weyl mixtes,
 Bol. Acad. Nac. Cienc. (C\'ordoba)
{ 65} (2000),
29--43.


\bibitem{AOVB}
J.
Alev, A. Ooms, M. Van den Bergh,
A class of counterexamples to the Gelfand-Kirillov conjecture,
Trans. Amer. Math. Soc. { 348} (1996), no. 5, 1709--1716. 








\bibitem{A}
V.A. Artamonov,
The skew field of rational quantum functions,
Russian Math. Surveys { 54} (1999), no. 4,
825--827.

\bibitem{A2}
$\soul\esp$,
General quantum polynomials: irreducible
    modules and Morita-equivalence,
    Izv. RAN, ser. math. { 63} (1999), no. 5, 3-36.



\bibitem{AVV}  M. Awami, M. Van den Bergh, F. Van Oystaeyen, 
 Note on derivations of graded rings and classification of differential polynomial rings, {\it in}
 Deuxième
  Contact Franco-Belge en Algèbre (Faulx-les-Tombes, 1987). Bull. Soc. Math. Belg. Sér. A 40 (1988), no. 2, 175--183. 










\bibitem{BG}
K.A. Brown, K.R. Goodearl,
Lectures on algebraic quantum groups.
Advanced Course in Math. CRM Barcelona, Vol. 2.
Birkh\"auser Verlag, Basel, 2002.

\bibitem{Cal}
  P. Caldero,
On the Gelfand-Kirillov conjecture for quantum algebras, Proc. Amer. Math. Soc. { 128} (2000), no. 4, 943--951. 

\bibitem{Ca1}
G. Cauchon,
 S\'eries formelles crois\'ees,
 J. Pure  Appl. Algebra { 107} (1996),
no. 2-3,
 153--169.


\bibitem{Ca2}
$\soul\esp$,
Effacement des d\'erivations et spectres premiers des
alg\`ebres quantiques,
 J. Algebra { 260} (2003), no. 2, 476--518. 




\bibitem{CGG}
L.J. Corwin, I.M. Gelfand, R. Goodman, 
Quadratic algebras and skew-fields,
in Representation theory and analysis on homogeneous spaces (New Brunswick, NJ, 1993).
Contemp. Math., Vol. 117.
Amer. Math. Soc., Providence, RI, 1994. pp 217--225.
 
 
\bibitem{De}
 E.E. Demidov,
{Some aspects of the theory of quantum groups},
Russian Math. Surveys
{ 38}
(1993), no. 6,
41--79.


 
 

 
\bibitem{FKK}
 H. Fujita, E.E. Kirkman, J. Kuzmanovich,
{Global and Krull dimension of quantum Weyl algebras},
J. Algebra
{ 216}
(1999), no. 2,
405--416.

 
 
\bibitem{F}
N. 
Fukuda, 
{Inverse and direct images for quantum Weyl algebras}, 
J. Algebra { 236} (2001), no. 2, 471--501. 

 
\bibitem{GK}
I.M. Gelfand, A.A. Kirillov,
{Sur les corps li\'es aux alg\`ebres enveloppantes
des alg\`ebres de Lie},
 Inst. Hautes Etudes Sci. Publ. Math.
{ 31}
 (1966), 509--523.
          
                                      
\bibitem{GZ}
A. Giaquinto, J.J. Zhang,
{Quantum Weyl algebras},
J. Algebra
{ 174}
(1995), no. 3,
861--881.





\bibitem{HK}
C. Huh, C.O. Kim,  
{Gelfand-Kirillov dimension of skew polynomial rings of automorphism type},
Comm. Algebra { 24} (1996), no. 7, 2317--2323. 



\bibitem{IM}
 K.  Iohara, F. Malikov, {Rings of skew polynomials and Gelfand-Kirillov conjecture for quantum groups}, Comm. Math. Phys. { 164} (1994), no. 2, 217--237. 


\bibitem{Jslqwa}
D.A. Jordan, 
 {A simple localization of the quantized Weyl algebra},
 J. Algebra { 174} (1995),
no. 1,
 267--281.

\bibitem{Jed}
$\soul{\esp }$,
{The graded algebra generated by two Eulerian derivatives},
Algebr. Represent. Theory { 4} (2001), no. 3, 249--275. 





\bibitem{Josep}
A. Joseph,
{Sur une conjecture de Feigin},
C. R. Acad. Sci. Paris Sér. I Math. { 320} (1995), no. 12, 1441--1444. 
 







\bibitem{MCP}
J.C. MacConnell, J.J. Pettit,
 {Crossed products and multiplicative analogues of Weyl algebras},
 J. London Math. Soc. (2) { 38} (1988),
no. 1,
 47--55.

\bibitem{Malt}
G.   Maltsiniotis, 
{Groupes quantiques et structures différentielles}, 
 C. R. Acad. Sci. Paris S\'er. I Math. { 311} (1990), no. 12, 831--834.


\bibitem{Pa1}
A.N. Panov, 
 {Skew fields of twisted rational functions and the skew field of rational functions on $GL_q(n,K)$},
St. Petersburg Math. J. { 7} (1996),
no. 1,
 129--143.



\bibitem{Pa2}
$\soul\esp$,
 {Fields of fractions of quantum solvable algebras},
 J. Algebra { 236} (2001),
no. 1, 
 110--121.


\bibitem{Retq}
L. Richard,
{Sur les endomorphismes des tores quantiques},
  Comm. Algebra. { 30} (2002), no. 11, 5281--5304.
  
\bibitem{these}
$\soul{\esp}$,
{Equivalence rationnelle et homologie de Hochschild pour certaines alg\`ebres polynomiales  classiques et quantiques}, th\`ese de doctorat de l'Universit\'e Blaise Pascal, 2002.

\bibitem{Rhhqw}
$\soul{\esp}$,
{Hochschild homology and cohomology of some classical and quantum noncommutative polynomial algebras}, \`a para\^\i tre dans J. Pure Appl. Algebra.


\bibitem{WZ}
J. Wess, B. Zumino,
{Covariant differential calculus on the quantum hyperplane},
Nuclear Phys. B Proc. Suppl. {\bf 18} (1990), 302--312.



\bibitem{Z}
 J.J. Zhang,
 {On Gelfand-Kirillov transcendence degree},
 Trans. Amer. Math. Soc. { 348} (1996), no. 7, 2867--2899. 


\end{thebibliography}
\end{document}